\documentclass[hidelinks,onefignum,onetabnum,nohypdvips]{siamart251216}

\usepackage[utf8]{inputenc}
\usepackage{amssymb,amsfonts}
\usepackage{mathrsfs}
\usepackage{bm}
\usepackage{graphicx,epstopdf}
\ifpdf
  \DeclareGraphicsExtensions{.eps,.pdf,.png,.jpg}
\else
  \DeclareGraphicsExtensions{.eps}
\fi
\graphicspath{{media/}{figures/}}
\usepackage[caption=false]{subfig}   
\let\subfigure\subfloat              
\usepackage{algpseudocode}           
\usepackage{tikz,pgfplots}
\pgfplotsset{compat=1.18}
\usepackage{braket}
\usepackage{booktabs}
\usepackage{multirow}
\usepackage{makecell}
\usepackage{cases}
\usepackage{calligra}
\usepackage{enumerate}
\usepackage[normalem]{ulem}
\usepackage{cite}

\newsiamremark{remark}{Remark}

\newcommand\bx{\boldsymbol{x}}
\newcommand\bn{\boldsymbol{n}}
\newcommand\bu{\boldsymbol{u}}
\newcommand\bv{{\boldsymbol{v}}}

\newcommand\bg{\boldsymbol{g}}

\newcommand\bh{\boldsymbol{h}}
\newcommand\bbR{\mathbb{R}}

\newcommand\bbS{\mathbb{S}}

\newcommand\bdeta{\boldsymbol{\eta}}

\newcommand\Kn{{\rm Kn}}

\newcommand\mQ{\mathcal{Q}}

\newcommand{\bz}{{\bm 0}}

\headers{A fast Laguerre--Fourier spectral method for polyatomic gases}{P. Zhang, Y. Wang, and Z. Cai}

\title{A fast Laguerre--Fourier spectral method for the polyatomic Boltzmann collision operator with continuous internal energy%
  \thanks{\funding{The work of Zhenning Cai was supported by the Academic Research Fund of the Ministry of Education of Singapore under grant A-8002392-00-00. The work of Yanli Wang was partially supported by the Science Challenge program under Grant No. TZ2026004. The work of Pei Zhang was supported by the China Scholarship Council (No. 202404890001).}}}

\author{Pei Zhang\thanks{Beijing Computational Science Research Center, Beijing, China
  (\email{zhangpei@csrc.ac.cn}, \email{ylwang@csrc.ac.cn}).}
\and Yanli Wang\footnotemark[2]
\and Zhenning Cai\thanks{Department of Mathematics, National University of Singapore, Singapore 119076
  (\email{matcz@nus.edu.sg}).}}

\ifpdf
\hypersetup{
  pdftitle={A fast Laguerre--Fourier spectral method for the polyatomic Boltzmann collision operator with continuous internal energy},
  pdfauthor={Pei Zhang, Yanli Wang, and Zhenning Cai}
}
\fi

\begin{document}
\maketitle

\begin{abstract}
The polyatomic Boltzmann equation extends the classical kinetic description to gases whose molecules possess internal energy in addition to their translational motion. The distribution function then has one more variable than in the monatomic case, and the collision operator becomes an eight-fold integral whose evaluation dominates the cost of solving the equation. In this work, we propose a fast Laguerre--Fourier spectral method for the polyatomic Boltzmann equation with continuous internal energy, which evaluates the collision operator in $O(K^2MN^4\log N + K^3MN^4)$ online operations after a solution-independent precomputation, where $N$ is the number of Fourier modes in each velocity direction, $K$ is the number of Laguerre modes and $M$ is the number of quadrature points on the sphere. Specifically, the distribution function is expanded in generalized Laguerre functions in the internal energy and in Fourier modes in velocity. A sequence of changes of variables then separates the internal-energy integrals from the velocity convolutions, so that the former no longer involve the distribution function and are computed once and stored. Numerical experiments are presented for both the spatially homogeneous and the spatially inhomogeneous Boltzmann equation. They confirm the spectral accuracy of the proposed method and its agreement with the direct simulation Monte Carlo method.
\end{abstract}




\section{Introduction}
\label{sec:introduction}

Rarefied gas dynamics has wide applications in astronautics, micro/nano flows and other fields \cite{shen2006rarefied}. The degree of rarefaction is characterized by the Knudsen number $\mathrm{Kn}$, defined as the ratio of the molecular mean free path to a characteristic length of the flow. When $\mathrm{Kn}\ll1$, the gas is well described by the continuum-based Euler and Navier--Stokes equations. As $\mathrm{Kn}$ increases, however, the continuum hypothesis breaks down, and the Boltzmann equation is required to describe the dynamics of the gas. 

The Boltzmann equation is an integro-differential equation that couples transport of the gas molecules to their binary collisions and governs the evolution of the velocity distribution function. The main difficulties in solving the Boltzmann equation are the high dimensionality of the distribution function and the complexity of the collision operator. The collision operator is a nonlinear integral operator that describes binary collisions between molecules. Whereas monatomic molecules exchange only translational energy in collisions, those of a polyatomic gas also carry internal rotational and vibrational energy that is exchanged during collisions, so that the collision operator additionally integrates over the internal energies and the energy-partition parameters. A number of studies have modeled this internal energy as a continuous variable and analyzed the resulting operator \cite{borgnakke1975statistical,bourgat1994microreversible,djordjic2021consistent,djordjic2023boltzmann}. Developing efficient numerical methods for this complex collision operator is therefore of both theoretical and practical importance.

Numerical methods for the collision operator mainly fall into two categories: stochastic methods and deterministic methods. The direct simulation Monte Carlo (DSMC) method \cite{bird1994molecular} is the most widely used stochastic particle method, representing the distribution by simulation particles and reproducing the collision dynamics through random binary collisions. It is highly efficient in the strongly rarefied regime and extends naturally to polyatomic gases \cite{borgnakke1975statistical,djordjic2023boltzmann}, but its statistical noise is expensive to control near the continuum limit and in low-speed or unsteady flows. Deterministic methods, including the discrete velocity method (DVM) \cite{broadwell1964study,goldstein1989investigations} and the spectral method \cite{pareschi1996fourier,pareschi2000numerical}, instead discretize the distribution function directly and are free of statistical noise. The DVM represents the distribution on a fixed velocity grid and computes the collision integral by numerical quadrature over that grid \cite{buet1996discrete}. Since the postcollision velocities rarely coincide with grid points, they must be interpolated back onto the grid, which makes the collision operator expensive to compute. The spectral method expands the distribution in a global basis and computes the collision term in transform space, which, combined with the fast Fourier transform, can make it considerably more efficient. Depending on the basis functions, spectral methods include the Fourier, Hermite \cite{wang2019approximation,li2023hermite}, Burnett \cite{cai2020burnett,hu2020burnett}, and mapped Chebyshev \cite{hu2020petrov,hu2022fast} spectral methods. Formulated on a truncated and periodized velocity domain, the Fourier spectral method can exploit the fast Fourier transform and has become the most widely used spectral approach for the collision operator.

For monatomic gases the Fourier spectral method has been developed extensively. The direct Fourier spectral method \cite{pareschi2000numerical} attains spectral accuracy but has a computational cost of $O(N^6)$, where $N$ is the number of modes in each velocity direction. This cost was reduced by the fast spectral methods \cite{mouhot2006fast,filbet2006solving,gamba2017fast}, which utilize the convolutional structure of the collision term together with a separable approximation of the collision kernel to bring the complexity down to $O(M N^4\log N)$, with $M$ the number of quadrature points on the sphere. Beyond the elastic monatomic operator \cite{wu2013deterministic}, the fast spectral framework has been extended to many other collision models, including inelastic and granular gases \cite{hu2019fast}, the quantum Uehling--Uhlenbeck operator \cite{hu2012fast}, the Fokker--Planck--Landau operator \cite{pareschi2000fokker}, gas mixtures \cite{jaiswal2019multispecies}, and the linearized collision operator \cite{yin2025fast}. 

The polyatomic collision operator is far more complex than its monatomic counterpart. Beyond the velocity variable, it integrates over the internal energy of the collision partner and the Borgnakke--Larsen energy-partition parameters. These additional integrals remain coupled to the velocity computation, so a direct quadrature evaluation of the operator is prohibitively expensive. To avoid this computational cost, deterministic simulations of polyatomic gases often adopt simplified kinetic models, such as the polyatomic ES-BGK model \cite{mathiaud2022esbgk,arima2024novel} and the Rykov model and its extensions \cite{rykov1975model,zeng2023general}. These models replace the full collision integral with simplified collision operators that account for energy exchange between translational and internal modes, thereby reducing the complexity of the numerical computation. Deterministic methods for the polyatomic collision operator have accordingly received little study. For discrete energy levels, the spectral-Lagrangian solver of \cite{munafo2014spectral} computes the operator by a direct Fourier spectral method at a cost of $O(K^4 N^6)$, where $K$ is the number of internal energy levels. For the continuous internal energy operator, no deterministic method has to our knowledge been developed. 

In this paper we develop such a method, a fast Laguerre--Fourier spectral method (FLFSM) for the polyatomic Boltzmann collision operator with continuous internal energy. Guided by the structure of the polyatomic equilibrium, we discretize the internal energy with generalized Laguerre basis functions \cite{shen2011spectral} and the velocity with Fourier basis functions, which yields a biorthogonal Laguerre--Fourier discretization of the collision operator. In the velocity variable this discretization retains the monatomic structure, so the monatomic fast spectral method, which exploits the convolution structure, still applies. However, the coupling between the internal energy and the velocity makes a direct application of this method expensive, since the integrals over the internal energies and the energy-partition parameters have to be evaluated online together with the convolutions. To reduce this cost, we derive a separable structure in which the internal energy and energy-partition integrals decouple from the velocity convolution and collapse into solution-independent coefficients computed once and stored offline, lowering the online complexity to $O(K^2 M N^4\log N+K^3MN^4)$, compared with $O(K^3N^6)$ for the direct spectral method. We further validate the method numerically. Because no exact solution of the polyatomic Boltzmann equation is known at the distribution level, we derive an analytic two-temperature relaxation law for Maxwell molecules in the spatially homogeneous case, which serves as an exact macroscopic reference. We then confirm the method's spectral convergence and compare it against DSMC on spatially homogeneous relaxation and on one-dimensional Fourier and Couette flows.

The remainder of the paper is organized as follows. Sec.~\ref{sec:polyatomic_fourier} introduces the polyatomic Boltzmann equation with continuous internal energy. Sec.~\ref{sec:fast_fourier} develops the direct and fast Laguerre--Fourier spectral methods, together with the algorithm and its complexity analysis. Sec.~\ref{sec:numerical} presents the numerical results, and Sec.~\ref{sec:conclusion} concludes.

\section{The polyatomic Boltzmann equation and its properties}
\label{sec:polyatomic_fourier}

\label{subsec:polyatomic_equation}

The polyatomic Boltzmann equation describes the kinetic evolution of a dilute gas whose molecules possess internal degrees of freedom in addition to their translational motion. The state of the gas is described by a distribution function $f(t,\bx,\bv,I)$, where $t$ is the time variable, $\bx$ is the spatial location, $\bv\in\bbR^3$ is the molecular velocity, and $I\in\bbR_+$ is a continuous internal energy variable that accounts for the internal degrees of freedom of the molecule, such as rotation \cite{bourgat1994microreversible,borgnakke1975statistical}. The polyatomic Boltzmann equation has the form
\begin{equation}
    \frac{\partial f}{\partial t} + \bv\cdot\nabla_{\bx} f = \mQ(f,f)(t,\bx,\bv,I), \qquad t\in\bbR_+,\quad \bx\in\Omega , \label{eq:bolt}
\end{equation}
where $\Omega\subset\bbR^3$ is the spatial domain. The right-hand side $\mQ(f,f)$ describes binary collisions in which the translational and internal energies may be redistributed \cite{djordjic2023boltzmann}, and takes the form
\begin{equation}
    \begin{aligned}
        \mQ(f,f)(\bv,I) = \int_{\bbR^3}\!\int_0^\infty\!\int_{\bbS^2}\!\int_0^1\!\int_0^1 & \left[ f(\bv',I')\,f(\bv_*',I_*') \left(\frac{I I_*}{I'I_*'}\right)^{\frac{\delta}{2}-1} - f(\bv,I)\,f(\bv_*,I_*) \right] \\
        &\times \mathcal B\, H_\delta(r,R)\, dR\,dr\,d\boldsymbol{\sigma}\,dI_*\,d\bv_* .
    \end{aligned}
    \label{eq:coll}
\end{equation}
Here, $\delta>0$ is the number of internal degrees of freedom, and a collision maps the precollision states $(\bv,I)$, $(\bv_*,I_*)$ to the postcollision states $(\bv',I')$, $(\bv_*',I_*')$ while conserving momentum and total energy,
\begin{equation}
    \bv+\bv_*=\bv'+\bv_*', \qquad \tfrac{1}{2}|\bv|^2+I+\tfrac{1}{2}|\bv_*|^2+I_* = \tfrac{1}{2}|\bv'|^2+I'+\tfrac{1}{2}|\bv_*'|^2+I_*' . \label{eq:conser}
\end{equation}
Let $\bh=(\bv+\bv_*)/2$ and $\bg=\bv-\bv_*$ be the mean and relative velocities. In these variables the conservation laws \eqref{eq:conser} become the invariance of $\bh$ and of the total energy
\begin{equation}
    E=\tfrac{1}{4}|\bg|^2+I+I_* . \label{eq:poly_total_energy}
\end{equation}
In the Borgnakke--Larsen parametrization \cite{borgnakke1975statistical}, the parameter $R\in[0,1]$ in \eqref{eq:coll} is an energy-partition variable and the postcollision relative kinetic energy is $\tfrac{1}{4}|\bg'|^2=RE$, from which the postcollision velocities are obtained as
\begin{equation}
    \bv'=\bh+\sqrt{R E}\,\boldsymbol{\sigma}, \qquad \bv_*'=\bh-\sqrt{R E}\,\boldsymbol{\sigma}, \label{eq:poly_post_velocity}
\end{equation}
where $\boldsymbol{\sigma}\in\bbS^2$ is a unit vector giving the direction of the postcollision relative velocity. The parameter $r\in[0,1]$ in \eqref{eq:coll} sets the proportion in which the total postcollision internal energy $(1-R)E$ is distributed between the two molecules, whose postcollision internal energies are then
\begin{equation}
    I'=r(1-R)E, \qquad I_*'=(1-r)(1-R)E . \label{eq:poly_post_internal_energy}
\end{equation}
The weight $H_\delta(r,R)$ in \eqref{eq:coll} has the form
\begin{equation}
    H_\delta(r,R) = [r(1-r)]^{\frac{\delta}{2}-1} (1-R)^{\delta-1} R^{1/2}, 
    \label{eq:Hdelta}
\end{equation}
which renders the Borgnakke--Larsen collision measure invariant under the exchange of pre- and postcollision variables \cite{djordjic2021consistent}.

The nonnegative kernel $\mathcal B=\mathcal B(\bg,I,I_*,r,R,\boldsymbol{\sigma})$ in \eqref{eq:coll} is determined by the intermolecular interaction. In this work we consider collision kernels of the form
\begin{equation}
    \mathcal{B}(|\bg|,R,\hat{\bg}\cdot\boldsymbol{\sigma}) = \mathfrak b(\hat{\bg}\cdot\boldsymbol{\sigma})\,|\bg|^{\zeta}\,R^{\zeta/2}, \qquad \zeta\geqslant 0, \label{eq:vhs_kernel}
\end{equation}
where $\hat{\bg}=\bg/|\bg|$, $\hat{\bg}\cdot\boldsymbol{\sigma}=\cos\chi$ with $\chi$ the scattering angle between the pre- and postcollision relative velocities, and $\mathfrak b\geqslant0$ is the angular part of the kernel. In particular, for $\mathfrak b\equiv B_0>0$, the kernel \eqref{eq:vhs_kernel} reduces to the variable-hard-sphere (VHS) type kernel for continuous internal energy \cite{djordjic2023boltzmann},
\begin{equation}
    \mathcal{B}(|\bg|,R) = B_0\,|\bg|^{\zeta}\,R^{\zeta/2}, \label{eq:vhs_iso}
\end{equation}
where $\zeta=0$ corresponds to Maxwell molecules and $\zeta=1$ to hard spheres.

As in the monatomic case, the collision operator \eqref{eq:coll} conserves mass, momentum, and total energy:
\begin{equation}
    \int_{\bbR^3}\int_0^\infty \mQ(f,f)
    \begin{pmatrix}
    1\\[1mm]
    \bv\\[1mm]
    \frac{1}{2}|\bv|^2+I
    \end{pmatrix}
    dI\,d\bv = \bz . \label{eq:poly_conservation}
\end{equation}
The density $\rho$, velocity $\bu$, and temperature $\theta$ are defined as the moments
\begin{equation}
\begin{gathered}
    \rho = \!\int_{\bbR^3}\!\int_0^\infty\! f\,dI\,d\bv, \qquad \rho\bu = \!\int_{\bbR^3}\!\int_0^\infty\! \bv f\,dI\,d\bv, \\
    \frac{3+\delta}{2}\rho\theta =\!\int_{\bbR^3}\!\int_0^\infty\!\! \Big(\frac{1}{2}|\bv-\bu|^2+I\Big) f\,dI\,d\bv .
\end{gathered}
    \label{eq:poly_moments}
\end{equation}
We also distinguish the translational and internal temperatures,
\begin{equation}
    \frac{3}{2}\rho\theta_{\mathrm{tr}} = \int_{\bbR^3}\int_0^\infty \frac{1}{2}|\bv-\bu|^2 f\,dI\,d\bv , \qquad \frac{\delta}{2}\rho\theta_{\mathrm{int}} = \int_{\bbR^3}\int_0^\infty I f\,dI\,d\bv , \label{eq:poly_temperatures}
\end{equation}
which coincide at equilibrium, $\theta_{\mathrm{tr}}=\theta_{\mathrm{int}}=\theta$; their relaxation toward the common value $\theta$ is derived in App.~\ref{subsec:two_temperature_law}. The stress tensor $P_{ij}$ and heat flux $q_i$ are higher-order moments of the distribution function, given by
\begin{equation}
\begin{aligned}
    P_{ij} &= \int_{\bbR^3}\!\int_0^\infty (v_i-u_i)(v_j-u_j)\, f\,dI\,d\bv , \\
    q_i &= \int_{\bbR^3}\!\int_0^\infty \Big(\tfrac12|\bv-\bu|^2+I\Big)(v_i-u_i)\, f\,dI\,d\bv , \qquad i,j=1,2,3 .
\end{aligned}
    \label{eq:poly_stress_heat}
\end{equation}
The equilibrium distribution is a Maxwellian in $\bv$ multiplied by a Gamma-type distribution in $I$,
\begin{equation}
    M(\bv,I) = \rho \frac{ I^{\frac{\delta}{2}-1} \exp\!\left( -\dfrac{\frac{1}{2}|\bv-\bu|^2+I}{\theta} \right) } {(2\pi\theta)^{3/2}\theta^{\delta/2}\Gamma(\delta/2)} . \label{eq:poly_maxwellian}
\end{equation}

\section{A fast Laguerre--Fourier spectral method}
\label{sec:fast_fourier}

In this section, we construct the fast Laguerre--Fourier spectral method for the collision operator \eqref{eq:coll}, where the internal energy is discretized by generalized Laguerre functions and the velocity by Fourier modes. Sec.~\ref{subsec:fourier_spectral_formulation} introduces the Laguerre--Fourier discretization and the direct spectral method, Sec.~\ref{subsec:fast_separable} develops the fast algorithm, and Sec.~\ref{subsec:fast_gain_loss} presents the fully discrete scheme and its complexity.

\subsection{Laguerre--Fourier discretization and the direct spectral method}
\label{subsec:fourier_spectral_formulation}

Motivated by the form of the equilibrium \eqref{eq:poly_maxwellian}, we approximate the dependence of $f$ on $I$ by generalized Laguerre functions. The generalized Laguerre polynomials $\{L_n^{(\alpha)}(z)\}_{n\geqslant 0}$ are orthogonal on $[0,\infty)$ with respect to the weight $z^{\alpha}e^{-z}$ \cite{shen2011spectral},
\begin{equation}
    \int_0^\infty L_m^{(\alpha)}(z)\,L_n^{(\alpha)}(z)\,z^{\alpha}e^{-z}\,dz = \frac{\Gamma(n+\alpha+1)}{n!}\,\delta_{mn}, \qquad \alpha>-1 . \label{eq:lag_orth}
\end{equation}
Comparing this weight with the internal-energy factor \eqref{eq:poly_maxwellian}, we take $\alpha=\tfrac{\delta}{2}-1$ together with the scaling $z=I/\beta$, where $\beta>0$ is a scaling parameter.

We define the test functions as the generalized Laguerre functions
\begin{equation}
    \Psi_b(I)=e^{-I/(2\beta)}\,L_b^{(\alpha)}\!\big(I/\beta\big), \label{eq:psi_b}
\end{equation}
and the trial functions as
\begin{equation}
    \Phi_b(I)=\frac{1}{Z_b} \Big(\frac{I}{\beta}\Big)^{\alpha} \Psi_b(I), \qquad Z_b=\beta\,\frac{\Gamma(b+\alpha+1)}{b!} .
     \label{eq:phi_b}
\end{equation}
By \eqref{eq:lag_orth}, the trial and test functions satisfy the biorthogonality relation $\int_0^\infty \Psi_a(I)\Phi_b(I)\,dI=\delta_{ab}$. We expand the distribution function as
\begin{equation}
    f(\bv,I) = \sum_{b=0}^{\infty} f_b(\bv)\,\Phi_b(I), \label{eq:f_exp}
\end{equation}
where $f_b(\bv)$ is the coefficient of the $b$-th mode. By the biorthogonality relation, the coefficients are given by
\begin{equation}
    f_b(\bv) = \int_0^\infty f(\bv,I)\,\Psi_b(I)\,dI . \label{eq:fb_coef}
\end{equation}

Each coefficient $f_b(\bv)$ is further represented through its Fourier transform in $\bv$ \cite{pareschi2000numerical,mouhot2006fast,gamba2017fast}. We denote by $\widehat f_b(\bdeta)$ the Fourier transform of $f_b(\bv)$, defined by
\begin{equation}
    \widehat f_b(\bdeta) = \int_{\bbR^3} f_b(\bv)\,e^{-i\bv\cdot\bdeta}\,d\bv , \qquad f_b(\bv) = \frac{1}{(2\pi)^3} \int_{\bbR^3} \widehat f_b(\bdeta)\,e^{i\bv\cdot\bdeta}\,d\bdeta . \label{eq:v_ft}
\end{equation}
Applying the same discretization to the collision operator, we define its Laguerre--Fourier coefficients by
\begin{equation}
    \widehat \mQ_a(\boldsymbol{\xi}) = \int_{\bbR^3}\int_0^\infty \mQ(f,f)(\bv,I)\,\Psi_a(I)\, e^{-i\bv\cdot\boldsymbol{\xi}}\,dI\,d\bv . \label{eq:Q_proj}
\end{equation}
Substituting the collision operator \eqref{eq:coll} with the kernel \eqref{eq:vhs_kernel}, the internal-energy expansion \eqref{eq:f_exp}, and the velocity Fourier representation \eqref{eq:v_ft} of each coefficient $f_b$ into \eqref{eq:Q_proj}, and separating the gain and loss contributions of $\mQ$, we have $\widehat \mQ_a=\widehat \mQ_a^+-\widehat \mQ_a^-$, where the gain term is
\begin{equation}
\begin{aligned}
    \widehat \mQ_a^+(\boldsymbol{\xi}) &= \frac{1}{(2\pi)^6}\sum_{b,d=0}^{\infty} \int_{\bbR^3}\!\int_{\bbR^3}\!\int_{\bbR^3}\!\int_{\bbR^3}\!\int_0^\infty\!\int_0^\infty\!\int_{\bbS^2}\!\int_0^1\!\int_0^1  \widehat f_b(\bdeta)\,\widehat f_d(\boldsymbol{\zeta})\,\Psi_a(I) \\
    &\times \Phi_b(I')\Phi_d(I_*') \left(\frac{I I_*}{I'I_*'}\right)^{\frac{\delta}{2}-1} \mathcal B(|\bg|,R,\hat{\bg}\cdot\boldsymbol{\sigma})\,H_\delta(r,R) \\
    &\times e^{\,i(\bv'\cdot\bdeta+\bv_*'\cdot\boldsymbol{\zeta}-\bv\cdot\boldsymbol{\xi})} \, d\bdeta\,d\boldsymbol{\zeta}\,d\bv\,d\bv_*\,dI\,dI_*\,d\boldsymbol{\sigma}\,dr\,dR ,
\end{aligned}
    \label{eq:gn_ft}
\end{equation}
and the loss term is
\begin{equation}
\begin{aligned}
    \widehat \mQ_a^-(\boldsymbol{\xi}) &= \frac{1}{(2\pi)^6}\sum_{b,d=0}^{\infty} \int_{\bbR^3}\!\int_{\bbR^3}\!\int_{\bbR^3}\!\int_{\bbR^3}\!\int_0^\infty\!\int_0^\infty\!\int_{\bbS^2}\!\int_0^1\!\int_0^1  \widehat f_b(\bdeta)\,\widehat f_d(\boldsymbol{\zeta})\,\Psi_a(I)\,\Phi_b(I)\Phi_d(I_*) \\
    &\times \mathcal B(|\bg|,R,\hat{\bg}\cdot\boldsymbol{\sigma})\,H_\delta(r,R)\, e^{\,i(\bv\cdot\bdeta+\bv_*\cdot\boldsymbol{\zeta}-\bv\cdot\boldsymbol{\xi})} \, d\bdeta\,d\boldsymbol{\zeta}\,d\bv\,d\bv_*\,dI\,dI_*\,d\boldsymbol{\sigma}\,dr\,dR .
\end{aligned}
    \label{eq:ls_ft}
\end{equation}
Here $\bdeta,\boldsymbol{\zeta}\in\bbR^3$ are the Fourier variables of $f_b$ and $f_d$. In \eqref{eq:gn_ft} and \eqref{eq:ls_ft}, the distribution function enters only through the factor $\widehat f_b(\bdeta)\,\widehat f_d(\boldsymbol{\zeta})$.

In the direct Fourier spectral method \cite{pareschi2000numerical}, all the remaining integrals are precomputed and stored as a collision kernel. To this end, we change variables to the mean and relative velocities $(\bh,\bg)$, with $d\bv\,d\bv_*=d\bh\,d\bg$. By \eqref{eq:poly_post_velocity}, the phases in \eqref{eq:gn_ft} and \eqref{eq:ls_ft} can be written as
\begin{equation}
\begin{aligned}
    \bv'\cdot\bdeta+\bv_*'\cdot\boldsymbol{\zeta}-\bv\cdot\boldsymbol{\xi}
    &= \bh\cdot(\bdeta+\boldsymbol{\zeta}-\boldsymbol{\xi}) + \sqrt{RE}\,\boldsymbol{\sigma}\cdot(\bdeta-\boldsymbol{\zeta}) - \tfrac{1}{2}\,\bg\cdot\boldsymbol{\xi} ,\\
    \bv\cdot\bdeta+\bv_*\cdot\boldsymbol{\zeta}-\bv\cdot\boldsymbol{\xi}
    &= \bh\cdot(\bdeta+\boldsymbol{\zeta}-\boldsymbol{\xi}) + \tfrac{1}{2}\,\bg\cdot(\bdeta-\boldsymbol{\zeta}-\boldsymbol{\xi}) .
\end{aligned}
    \label{eq:phase_direct}
\end{equation}
The integral with respect to $\bh$ can be evaluated explicitly,
\begin{equation}
    \int_{\bbR^3}e^{i\bh\cdot(\bdeta+\boldsymbol{\zeta}-\boldsymbol{\xi})}\,d\bh=(2\pi)^3\delta_{\mathrm D}(\bdeta+\boldsymbol{\zeta}-\boldsymbol{\xi}),
    \label{eq:h_int}
\end{equation}
where $\delta_{\mathrm D}$ is the Dirac delta function. After integration in $\boldsymbol{\zeta}$, this enforces $\boldsymbol{\zeta}=\boldsymbol{\xi}-\bdeta$, and the gain term \eqref{eq:gn_ft} becomes
\begin{equation}
\begin{aligned}
    \widehat \mQ_a^+(\boldsymbol{\xi}) &= \frac{1}{(2\pi)^3}\sum_{b,d=0}^{\infty} \int_{\bbR^3}\!\int_{\bbR^3}\!\int_0^\infty\!\int_0^\infty\!\int_{\bbS^2}\!\int_0^1\!\int_0^1 \widehat f_b(\bdeta)\,\widehat f_d(\boldsymbol{\xi}-\bdeta)\,\Psi_a(I) \\
    &\times \Phi_b(I')\Phi_d(I_*') \left(\frac{I I_*}{I'I_*'}\right)^{\frac{\delta}{2}-1} \mathcal B(|\bg|,R,\hat{\bg}\cdot\boldsymbol{\sigma})\,H_\delta(r,R) \\
    &\times e^{\,i\sqrt{RE}\,\boldsymbol{\sigma}\cdot(2\bdeta-\boldsymbol{\xi})-\frac{i}{2}\bg\cdot\boldsymbol{\xi}}\, dR\,dr\,d\boldsymbol{\sigma}\,dI_*\,dI\,d\bg\,d\bdeta ,
\end{aligned}
    \label{eq:gn_hint}
\end{equation}
and the loss term \eqref{eq:ls_ft} becomes
\begin{equation}
\begin{aligned}
    \widehat \mQ_a^-(\boldsymbol{\xi}) &= \frac{1}{(2\pi)^3}\sum_{b,d=0}^{\infty} \int_{\bbR^3}\!\int_{\bbR^3}\!\int_0^\infty\!\int_0^\infty\!\int_{\bbS^2}\!\int_0^1\!\int_0^1 \widehat f_b(\bdeta)\,\widehat f_d(\boldsymbol{\xi}-\bdeta)\,\Psi_a(I)\,\Phi_b(I)\Phi_d(I_*) \\
    &\times \mathcal B(|\bg|,R,\hat{\bg}\cdot\boldsymbol{\sigma})\,H_\delta(r,R)\, e^{\,i\bg\cdot(\bdeta-\boldsymbol{\xi})}\, dR\,dr\,d\boldsymbol{\sigma}\,dI_*\,dI\,d\bg\,d\bdeta ,
\end{aligned}
    \label{eq:ls_hint}
\end{equation}
where $E$ is given by \eqref{eq:poly_total_energy}, and $I'$, $I_*'$ are given by \eqref{eq:poly_post_internal_energy}.

The integrals with respect to $(\bg,I,I_*,\boldsymbol{\sigma},r,R)$ in \eqref{eq:gn_hint} and \eqref{eq:ls_hint} are independent of $f$. Collecting them into a kernel, we obtain
\begin{equation}
    \widehat \mQ_a(\boldsymbol{\xi}) = \sum_{b,d=0}^{\infty} \int_{\bbR^3} \mathcal T_{a,b,d}(\bdeta,\boldsymbol{\xi})\, \widehat f_b(\bdeta)\,\widehat f_d(\boldsymbol{\xi}-\bdeta)\, d\bdeta ,
    \label{eq:Q_direct}
\end{equation}
where $\mathcal T_{a,b,d}=\mathcal T^+_{a,b,d}-\mathcal T^-_{a,b,d}$ with
\begin{equation}
\begin{aligned}
    \mathcal T^+_{a,b,d}(\bdeta,\boldsymbol{\xi}) &= \frac{1}{(2\pi)^3} \int_{\bbR^3}\!\int_0^\infty\!\int_0^\infty\!\int_{\bbS^2}\!\int_0^1\!\int_0^1 \Psi_a(I)\,\Phi_b(I')\Phi_d(I_*') \left(\frac{I I_*}{I'I_*'}\right)^{\frac{\delta}{2}-1} \\
    &\qquad\times \mathcal B(|\bg|,R,\hat{\bg}\cdot\boldsymbol{\sigma})\,H_\delta(r,R)\, e^{\,i\sqrt{RE}\,\boldsymbol{\sigma}\cdot(2\bdeta-\boldsymbol{\xi})-\frac{i}{2}\bg\cdot\boldsymbol{\xi}}\, dR\,dr\,d\boldsymbol{\sigma}\,dI_*\,dI\,d\bg ,\\
    \mathcal T^-_{a,b,d}(\bdeta,\boldsymbol{\xi}) &= \frac{1}{(2\pi)^3} \int_{\bbR^3}\!\int_0^\infty\!\int_0^\infty\!\int_{\bbS^2}\!\int_0^1\!\int_0^1 \Psi_a(I)\,\Phi_b(I)\Phi_d(I_*) \\
    &\qquad\times \mathcal B(|\bg|,R,\hat{\bg}\cdot\boldsymbol{\sigma})\,H_\delta(r,R)\, e^{\,i\bg\cdot(\bdeta-\boldsymbol{\xi})}\, dR\,dr\,d\boldsymbol{\sigma}\,dI_*\,dI\,d\bg .
\end{aligned}
    \label{eq:beta_direct}
\end{equation}
Since $\mathcal T_{a,b,d}$ is independent of $f$, it can be computed once in a preprocessing step and stored. Truncating the Laguerre expansion \eqref{eq:f_exp} to $K$ modes and using $N$ Fourier modes in each velocity direction, the evaluation of \eqref{eq:Q_direct} for all $a$ and $\boldsymbol{\xi}$ requires
\begin{equation}
    C_{\mathrm{dir}} = O\!\left( K^3 N^6 \right)
    \label{eq:cost_direct}
\end{equation}
operations, and storing $\mathcal T_{a,b,d}$ also requires $O(K^3N^6)$ memory. As $\mathcal T_{a,b,d}(\bdeta,\boldsymbol{\xi})$ is a nine-dimensional tensor, its storage exceeds the memory of standard workstations even for moderate values of $K$ and $N$. This motivates the development of a more efficient algorithm, which is the subject of the next subsection.

\subsection{The fast Laguerre--Fourier spectral method}
\label{subsec:fast_separable}

For the monatomic Boltzmann equation, the fast spectral method of \cite{gamba2017fast} precomputes the integral over the direction of the relative velocity and writes the collision term as a sum of convolutions, which reduces the cost from $O(N^6)$ to $O(MN^4\log N)$. We follow the same idea and first write \eqref{eq:gn_hint} and \eqref{eq:ls_hint} in convolution form.

\subsubsection{Convolution form of the collision term}
\label{subsubsec:conv_form}

We begin with the loss term \eqref{eq:ls_hint}. Its integrand depends on $I$ only through $\Psi_a(I)\Phi_b(I)$, so that the biorthogonality relation reduces the sum over $b$ to $b=a$. Moreover, a direct computation using \eqref{eq:phi_b} gives
\begin{equation}
    \int_0^\infty \Phi_d(I_*)\,dI_* = (-1)^d\,2^{\alpha+1} .
    \label{eq:Phi_int}
\end{equation}
The loss term \eqref{eq:ls_hint} therefore becomes
\begin{equation}
    \widehat \mQ_a^-(\boldsymbol{\xi})
    = \frac{1}{(2\pi)^3}\sum_{d=0}^{\infty}(-1)^d\,2^{\alpha+1}\int_{\bbR^3}
    \widehat f_a(\bdeta)\,\widehat f_d(\boldsymbol{\xi}-\bdeta)\,
    \widehat{\mathcal B}(\bdeta-\boldsymbol{\xi})\,d\bdeta ,
    \label{eq:ls_fin}
\end{equation}
where
\begin{equation}
    \widehat{\mathcal B}(\bdeta-\boldsymbol{\xi}) = \int_{\bbR^3}\!\int_{\bbS^2}\!\int_0^1\!\int_0^1 \mathcal B(|\bg|,R,\hat{\bg}\cdot\boldsymbol{\sigma})\,H_\delta(r,R)\, e^{i\bg\cdot(\bdeta-\boldsymbol{\xi})}\, dR\,dr\,d\boldsymbol{\sigma}\,d\bg .
    \label{eq:ls_ker}
\end{equation}

We now turn to the gain term \eqref{eq:gn_hint}. We write the relative velocity in spherical coordinates, $\bg=g\boldsymbol{\omega}$ with $g=|\bg|$ and $\boldsymbol{\omega}\in\bbS^2$, so that $d\bg=g^2\,dg\,d\boldsymbol{\omega}$. Since $E=\tfrac14 g^2+I+I_*$ depends only on $g$, the direction $\boldsymbol{\omega}$ enters the integrand of \eqref{eq:gn_hint} only through $\mathcal B(g,R,\boldsymbol{\omega}\cdot\boldsymbol{\sigma})\,e^{-\frac{i}{2}g\boldsymbol{\omega}\cdot\boldsymbol{\xi}}$. As in \cite{gamba2017fast}, the integral with respect to $\boldsymbol{\omega}$ is collected into the weight function
\begin{equation}
    F(g,R,\boldsymbol{\sigma},\boldsymbol{\xi}) = g^2\int_{\bbS^2} \mathcal B(g,R,\boldsymbol{\omega}\cdot\boldsymbol{\sigma})\, e^{-\frac{i}{2}g\boldsymbol{\omega}\cdot\boldsymbol{\xi}}\, d\boldsymbol{\omega} ,
    \label{eq:F_weight}
\end{equation}
which is independent of the distribution function. For the VHS-type kernel \eqref{eq:vhs_iso}, $F$ has the explicit form
\begin{equation}
    F(g,R,\boldsymbol{\xi}) = 4\pi B_0\,g^{2+\zeta}R^{\zeta/2}\operatorname{sinc}\!\Big(\frac{g|\boldsymbol{\xi}|}{2}\Big), \qquad \operatorname{sinc}(x)=\frac{\sin x}{x} .
    \label{eq:F_vhs}
\end{equation}
The gain term \eqref{eq:gn_hint} then becomes
\begin{equation}
\begin{aligned}
    \widehat \mQ_a^+(\boldsymbol{\xi}) &= \frac{1}{(2\pi)^3} \sum_{b,d=0}^{\infty} \int_0^\infty\!\int_{\bbS^2}\! \int_0^\infty\!\int_0^\infty\!\int_0^1\!\int_0^1 F(g,R,\boldsymbol{\sigma},\boldsymbol{\xi})\, \Psi_a(I)\,\Phi_b(I')\Phi_d(I_*') \\
    &\times \left(\frac{I I_*}{I'I_*'}\right)^{\frac{\delta}{2}-1} H_\delta(r,R)\, \mathcal C_{b,d}\big(2\sqrt{RE},\boldsymbol{\sigma},\boldsymbol{\xi}\big)\, dR\,dr\,dI_*\,dI\,d\boldsymbol{\sigma}\,dg ,
\end{aligned}
    \label{eq:gn_F}
\end{equation}
where
\begin{equation}
    \mathcal C_{b,d}(\lambda,\boldsymbol{\sigma},\boldsymbol{\xi}) = \int_{\bbR^3} \big[\widehat f_b(\bdeta)\, e^{\frac{i}{2}\lambda\boldsymbol{\sigma}\cdot\bdeta}\big]\, \big[\widehat f_d(\boldsymbol{\xi}-\bdeta)\, e^{-\frac{i}{2}\lambda\boldsymbol{\sigma}\cdot(\boldsymbol{\xi}-\bdeta)}\big]\,d\bdeta .
    \label{eq:C_conv}
\end{equation}
For fixed $b$, $d$, $\lambda$, and $\boldsymbol{\sigma}$, \eqref{eq:C_conv} is the convolution of $\widehat f_b\,e^{\frac{i}{2}\lambda\boldsymbol{\sigma}\cdot\bdeta}$ and $\widehat f_d\,e^{-\frac{i}{2}\lambda\boldsymbol{\sigma}\cdot\bdeta}$ evaluated at $\boldsymbol{\xi}$.

Finally, the integral with respect to $r$ in \eqref{eq:gn_F} can be carried out in advance. By \eqref{eq:poly_post_internal_energy} and \eqref{eq:Hdelta},
\begin{equation}
    \left(\frac{I I_*}{I'I_*'}\right)^{\frac{\delta}{2}-1} H_\delta(r,R) = (II_*)^{\frac{\delta}{2}-1}(1-R)R^{1/2}E^{2-\delta},
\end{equation}
which is independent of $r$, so that $r$ enters \eqref{eq:gn_F} only through $\Phi_b(I')\Phi_d(I_*')$. We thus obtain
\begin{equation}
\begin{aligned}
    \widehat \mQ_a^+(\boldsymbol{\xi}) &= \frac{1}{(2\pi)^3} \sum_{b,d=0}^{\infty} \int_0^\infty\!\int_{\bbS^2}\! \int_0^\infty\!\int_0^\infty\!\int_0^1 F(g,R,\boldsymbol{\sigma},\boldsymbol{\xi})\, \Psi_a(I)\,(II_*)^{\frac{\delta}{2}-1} \\
    &\times (1-R)R^{1/2}E^{2-\delta}\, G_{b,d}\big((1-R)E\big)\, \mathcal C_{b,d}\big(2\sqrt{RE},\boldsymbol{\sigma},\boldsymbol{\xi}\big)\, dR\,dI_*\,dI\,d\boldsymbol{\sigma}\,dg ,
\end{aligned}
    \label{eq:gn_r}
\end{equation}
where
\begin{equation}
    G_{b,d}(Y) = \int_0^1 \Phi_b(rY)\,\Phi_d\big((1-r)Y\big)\,dr . \label{eq:G_bd}
\end{equation}

We now analyze the computational cost of evaluating \eqref{eq:ls_fin} and \eqref{eq:gn_r}, where the Laguerre expansion is truncated to $K$ modes and $N$ Fourier modes are used in each velocity direction. In the loss term \eqref{eq:ls_fin}, the sum over $d$ can be carried out first, so that only one convolution is needed for each $a$, and the cost is $O(KN^3\log N)$.

The gain term \eqref{eq:gn_r} is evaluated in two steps. First, the convolutions $\mathcal C_{b,d}$ depend on the integration variables $(g,I,I_*,R)$ only through $\lambda=2\sqrt{RE}$. Discretizing $\lambda$ with $N_\lambda$ points and $\boldsymbol{\sigma}$ with $M$ points, all the convolutions are computed by fast Fourier transforms in $O(K^2MN_\lambda N^3\log N)$ operations. Second, the outer integrals are approximated by quadrature rules, with $N_g$ points for $g$, $M$ points for $\boldsymbol{\sigma}$, $N_R$ points for $R$, and $O(K)$ points for each of $I$ and $I_*$. At each of the $O(K^2MN_gN_R)$ quadrature nodes, the sum over $(b,d)$ costs $O(K^2N^3)$ operations for all $\boldsymbol{\xi}$, which dominates the subsequent summation for all $a$. Hence, the cost of the gain term is
\begin{equation}
    C_{\mathrm{conv}} = O\big(K^2MN_\lambda N^3\log N + K^4MN_gN_RN^3\big).
    \label{eq:cost_dir}
\end{equation}
In \eqref{eq:cost_dir}, the cost of the outer integrals is typically much larger than that of the convolutions, since $K^2N_gN_R\gg N_\lambda\log N$ in practice, and therefore dominates the total cost of the gain term. Moreover, $2\sqrt{RE}$ at the quadrature nodes generally does not coincide with the nodes of $\lambda$, so that $\mathcal C_{b,d}$ has to be interpolated. Both difficulties are removed by the changes of variables introduced below.

\subsubsection{Separation of the internal-energy integrals}
\label{subsubsec:separation}

In what follows, we introduce two changes of variables in \eqref{eq:gn_r}, after which the convolution $\mathcal C_{b,d}$ depends only on the integration variables $\lambda$ and $\boldsymbol{\sigma}$, and all the other integrals are moved into precomputed weights.

We first consider the integrals with respect to $(I,I_*)\in(0,\infty)^2$. In \eqref{eq:gn_r}, the convolution $\mathcal C_{b,d}$ depends on $I$ and $I_*$ only through the total energy $E=\tfrac14 g^2+I+I_*$, that is, only through the sum $I+I_*$. This suggests replacing $(I,I_*)$, which describe the precollision internal energies of the two molecules, by their total $U$ and the fraction $\tau$ of $U$ carried by the first molecule,
\begin{equation}
    U=I+I_*, \qquad \tau=\frac{I}{I+I_*}, \qquad dI\,dI_*=U\,d\tau\,dU .
    \label{eq:cv_Ut}
\end{equation}
Since $I=\tau U$, $I_*=(1-\tau)U$, and $E=\tfrac14 g^2+U$, substituting \eqref{eq:cv_Ut} into \eqref{eq:gn_r} gives
\begin{equation}
\begin{aligned}
    \widehat \mQ_a^+(\boldsymbol{\xi}) &= \frac{1}{(2\pi)^3} \sum_{b,d=0}^{\infty} \int_0^\infty\!\int_{\bbS^2}\! \int_0^\infty\!\int_0^1\!\int_0^1 F(g,R,\boldsymbol{\sigma},\boldsymbol{\xi})\, \Psi_a(\tau U)\,U^{\delta-1}[\tau(1-\tau)]^{\frac{\delta}{2}-1} \\
    &\times (1-R)R^{1/2}E^{2-\delta}\, G_{b,d}\big((1-R)E\big)\, \mathcal C_{b,d}\big(2\sqrt{RE},\boldsymbol{\sigma},\boldsymbol{\xi}\big)\, dR\,d\tau\,dU\,d\boldsymbol{\sigma}\,dg .
\end{aligned}
    \label{eq:gn_Ut}
\end{equation}
In \eqref{eq:gn_Ut}, $\tau$ appears only in $\Psi_a(\tau U)\,[\tau(1-\tau)]^{\frac{\delta}{2}-1}$, and the integral with respect to $\tau$ can be precomputed. The gain term then becomes
\begin{equation}
\begin{aligned}
    \widehat \mQ_a^+(\boldsymbol{\xi}) &= \frac{1}{(2\pi)^3} \sum_{b,d=0}^{\infty} \int_0^\infty\!\int_{\bbS^2}\! \int_0^\infty\!\int_0^1 F(g,R,\boldsymbol{\sigma},\boldsymbol{\xi})\, \mathcal P_a(U)\,(1-R)R^{1/2}E^{2-\delta} \\
    &\times G_{b,d}\big((1-R)E\big)\, \mathcal C_{b,d}\big(2\sqrt{RE},\boldsymbol{\sigma},\boldsymbol{\xi}\big)\, dR\,dU\,d\boldsymbol{\sigma}\,dg ,
\end{aligned}
    \label{eq:gn_U}
\end{equation}
where
\begin{equation}
    \mathcal P_a(U) = U^{\delta-1}\int_0^1 \Psi_a(\tau U)\,[\tau(1-\tau)]^{\frac{\delta}{2}-1}\,d\tau .
    \label{eq:P_a}
\end{equation}

In \eqref{eq:gn_U}, the convolution $\mathcal C_{b,d}$ still depends on $g$, $U$, and $R$ through $2\sqrt{RE}$. For fixed $g$, the variables $(U,R)$ describe how the energy is distributed in a collision: $\tfrac14 g^2$ and $U$ are the precollision relative translational and internal energies, whose sum is $E$, and $R$ determines how $E$ is split after the collision into the postcollision relative translational energy $RE$ and internal energy $(1-R)E$. The same partition can equally be described by the postcollision energies themselves. We therefore replace $(U,R)$ by
\begin{equation}
    \lambda=2\sqrt{RE}, \qquad Y=(1-R)E, \qquad E=\tfrac14 g^2+U ,
    \label{eq:cv_lY}
\end{equation}
where $\lambda=|\bg'|$ is the postcollision relative speed and $Y=I'+I_*'$ is the total postcollision internal energy. In the new variables, the convolution $\mathcal C_{b,d}$ depends only on $\lambda$, and the redistribution coefficient $G_{b,d}$ only on $Y$. The inverse transformation and the Jacobian are
\begin{equation}
    U=\frac{\lambda^2}{4}+Y-\frac{g^2}{4}, \qquad
    R=\frac{\lambda^2}{\lambda^2+4Y}, \qquad
    dU\,dR=\frac{2\lambda}{\lambda^2+4Y}\,d\lambda\,dY ,
    \label{eq:EUR}
\end{equation}
and the domain $U\geqslant0$, $R\in[0,1]$ becomes $\lambda\geqslant0$, $Y\geqslant Y_0(g,\lambda)$ with $Y_0(g,\lambda)=\max\big(0,\tfrac14 (g^2-\lambda^2)\big)$. Substituting \eqref{eq:EUR} into \eqref{eq:gn_U}, we obtain
\begin{equation}
\begin{aligned}
    \widehat \mQ_a^+(\boldsymbol{\xi}) &= \frac{1}{4(2\pi)^3} \sum_{b,d=0}^{\infty} \int_0^\infty\!\int_{\bbS^2}\!\int_0^\infty\!\int_{Y_0(g,\lambda)}^\infty \lambda^2\,Y\left(\frac{\lambda^2}{4}+Y\right)^{-\delta-\frac12} \\
    &\times F\Big(g,\frac{\lambda^2}{\lambda^2+4Y},\boldsymbol{\sigma},\boldsymbol{\xi}\Big)\, \mathcal P_a\Big(\frac{\lambda^2}{4}+Y-\frac{g^2}{4}\Big)\, G_{b,d}(Y)\, \mathcal C_{b,d}(\lambda,\boldsymbol{\sigma},\boldsymbol{\xi})\, dY\,d\lambda\,d\boldsymbol{\sigma}\,dg .
\end{aligned}
    \label{eq:gn_lY}
\end{equation}
Now $\mathcal C_{b,d}(\lambda,\boldsymbol{\sigma},\boldsymbol{\xi})$ is independent of $g$ and $Y$, and the integrals with respect to $g$ and $Y$ can be precomputed. The gain term thus takes the form
\begin{equation}
    \widehat \mQ_a^+(\boldsymbol{\xi}) = \frac{1}{(2\pi)^3} \sum_{b,d=0}^{\infty} \int_{\bbS^2}\!\int_0^\infty W_{a,b,d}(\lambda,\boldsymbol{\sigma},\boldsymbol{\xi})\, \mathcal C_{b,d}(\lambda,\boldsymbol{\sigma},\boldsymbol{\xi})\, d\lambda\,d\boldsymbol{\sigma} ,
    \label{eq:gn_fin}
\end{equation}
with the weight
\begin{equation}
    W_{a,b,d}(\lambda,\boldsymbol{\sigma},\boldsymbol{\xi}) = \frac{\lambda^2}{4} \int_0^\infty \mathcal{K}_{a,b,d}(g,\lambda,\boldsymbol{\sigma},\boldsymbol{\xi})\,dg ,
    \label{eq:W_abd}
\end{equation}
where
\begin{equation}
\begin{aligned}
    \mathcal{K}_{a,b,d}(g,\lambda,\boldsymbol{\sigma},\boldsymbol{\xi}) = \int_{Y_0(g,\lambda)}^\infty & Y\left(\frac{\lambda^2}{4}+Y\right)^{-\delta-\frac12} F\Big(g,\frac{\lambda^2}{\lambda^2+4Y},\boldsymbol{\sigma},\boldsymbol{\xi}\Big) \\
    &\times \mathcal P_a\Big(\frac{\lambda^2}{4}+Y-\frac{g^2}{4}\Big)\,G_{b,d}(Y)\,dY .
\end{aligned}
    \label{eq:K_abd}
\end{equation}
The weight $W_{a,b,d}$ is independent of the distribution function. For the VHS-type kernel \eqref{eq:vhs_iso}, it follows from \eqref{eq:F_vhs} that $W_{a,b,d}$ is independent of $\boldsymbol{\sigma}$. Apart from the sum over $(b,d)$, \eqref{eq:gn_fin} has the same structure as the gain term in the monatomic fast spectral method \cite{gamba2017fast}.

We now analyze the computational cost of evaluating \eqref{eq:gn_fin}. As in Sec.~\ref{subsubsec:conv_form}, the convolutions $\mathcal C_{b,d}$ are computed in $O(K^2MN_\lambda N^3\log N)$ operations. Since $\lambda$ is now an integration variable, the convolutions are only needed at the quadrature nodes in $\lambda$, and no interpolation is required. The remaining sums over $(b,d)$ and the quadrature nodes in $\lambda$ and $\boldsymbol{\sigma}$ for all $a$ and $\boldsymbol{\xi}$ cost $O(K^3MN_\lambda N^3)$ operations. For the VHS-type kernel, $W_{a,b,d}$ is independent of $\boldsymbol{\sigma}$, so that the convolutions can be summed over the quadrature nodes in $\boldsymbol{\sigma}$ first, and this cost reduces to $O(K^2MN_\lambda N^3+K^3N_\lambda N^3)$. Hence, the cost of the gain term is
\begin{equation}
    C_{\mathrm{fast}}^+ = O\big(K^2MN_\lambda N^3\log N + K^3MN_\lambda N^3\big),
    \label{eq:cost_sep}
\end{equation}
which no longer involves the quadratures in $g$, $R$, $I$, and $I_*$. Compared with \eqref{eq:cost_dir}, these quadratures have been moved entirely into the offline computation of $W_{a,b,d}$.

\subsection{The fully discrete scheme for the polyatomic collision term}
\label{subsec:fast_gain_loss}

Secs.~\ref{subsec:fourier_spectral_formulation} and~\ref{subsec:fast_separable} construct the continuous form of the fast Laguerre--Fourier spectral method for the polyatomic collision term. In this subsection we present the corresponding numerical method and analyze its complexity. Combining the reduced gain term \eqref{eq:gn_fin} and loss term \eqref{eq:ls_fin}, the projected collision coefficient is
\begin{equation}
\begin{aligned}
    \widehat \mQ_a(\boldsymbol{\xi}) = {}& \frac{1}{(2\pi)^3} \sum_{b,d=0}^{K-1} \int_{\bbS^2}\!\int_0^\infty W_{a,b,d}(\lambda,\boldsymbol{\sigma},\boldsymbol{\xi})\, \mathcal C_{b,d}(\lambda,\boldsymbol{\sigma},\boldsymbol{\xi})\, d\lambda\,d\boldsymbol{\sigma} \\
    &- \frac{1}{(2\pi)^3}\sum_{d=0}^{K-1}(-1)^d\,2^{\alpha+1}\int_{\bbR^3} \widehat f_a(\bdeta)\,\widehat f_d(\boldsymbol{\xi}-\bdeta)\, \widehat{\mathcal B}(\bdeta-\boldsymbol{\xi})\,d\bdeta .
\end{aligned}
    \label{eq:Q_cont}
\end{equation}

As in the monatomic Fourier spectral method, we truncate and periodize in velocity on the box
\begin{equation}
    \mathcal D_L=[-L,L]^3 .
\end{equation}
Let
\begin{equation}
    \mathbb K_N=\{\boldsymbol{k}=(k_1,k_2,k_3): -\frac{N}{2}\leqslant k_j\leqslant \frac{N}{2}\}, \qquad
    \boldsymbol{\xi}_{\boldsymbol{k}}=\frac{\pi}{L}\boldsymbol{k},
\end{equation}
where $N$ is even. For each Laguerre mode we use the Fourier representation
\begin{equation}
    f_b^N(\bv)=\frac{1}{(2L)^3}\sum_{\boldsymbol{k}\in\mathbb K_N}
    \widehat f_{b,\boldsymbol{k}}\,e^{i\boldsymbol{\xi}_{\boldsymbol{k}}\cdot\bv},
    \qquad
    \widehat f_{b,\boldsymbol{k}}=\frac{1}{c_{\boldsymbol{k}}}\Big(\frac{2L}{N}\Big)^3\sum_{\boldsymbol{j}\in\mathbb K_N}\frac{1}{c_{\boldsymbol{j}}}\,
    f_b(\bv_{\boldsymbol{j}})\,e^{-i\boldsymbol{\xi}_{\boldsymbol{k}}\cdot\bv_{\boldsymbol{j}}} ,
    \label{eq:v_fs}
\end{equation}
where $\bv_{\boldsymbol{j}}=\tfrac{2L}{N}\boldsymbol{j}$, $\boldsymbol{j}\in\mathbb K_N$, are the nodes of the uniform grid on $\mathcal D_L$, the nodes with $j_l=\pm N/2$ being identified by periodicity, $c_{\boldsymbol{k}}=c_{k_1}c_{k_2}c_{k_3}$, and $c_{k}=2$ if $k=\pm N/2$ and $c_k=1$ otherwise. Since the discrete Fourier coefficients are $N$-periodic in each index, the modes $k_j=\pm N/2$ share the same coefficient, and the factor $1/c_{\boldsymbol{k}}$ splits it equally between them. We denote by $\mathcal F_N$ the discrete Fourier transform in \eqref{eq:v_fs}, which maps the $(N+1)^3$ grid values $\{u(\bv_{\boldsymbol{j}})\}_{\boldsymbol{j}\in\mathbb K_N}$ of a periodic function $u$ to its $(N+1)^3$ coefficients $\{\widehat u_{\boldsymbol{k}}\}_{\boldsymbol{k}\in\mathbb K_N}$, and by $\mathcal F_N^{-1}$ its inverse, which evaluates the first relation in \eqref{eq:v_fs} at the grid nodes. Both are implemented by fast Fourier transforms.

In the computation, both the precollision relative speed $g$ and the postcollision relative speed $\lambda$ in \eqref{eq:gn_fin}--\eqref{eq:K_abd} have to be truncated. Specifically, we set $\mathcal B(g,\cdot)=0$ for $g>S$ and restrict $\lambda$ to $[0,\Lambda]$. In the monatomic case $g=\lambda$, and hence $\Lambda=S$. For polyatomic gases, translational and internal energy are exchanged in a collision, and \eqref{eq:EUR} gives
\begin{equation}
    g^2=\lambda^2+4(Y-U).
    \label{eq:g_lambda}
\end{equation}
Both speeds are nevertheless controlled by the total energy, which is conserved in a collision,
\begin{equation}
    E=\tfrac14 g^2+U=\tfrac14\lambda^2+Y .
    \label{eq:E_both}
\end{equation}
We next show that the choice $\Lambda=S$ also leads to a negligible truncation error in the polyatomic case. If $\lambda>S$, then $E>\tfrac14 S^2$ by \eqref{eq:E_both}, since $Y\geqslant0$. Assume that $f\leqslant C\,M$ for a constant $C>0$ and a Maxwellian $M$ of the form \eqref{eq:poly_maxwellian}. Then the factor $f(\bv',I')f(\bv_*',I_*')$ in the gain term is bounded by $C^2M(\bv',I')M(\bv_*',I_*')$. Since $\tfrac12|\bv'-\bu|^2+\tfrac12|\bv_*'-\bu|^2=|\bh-\bu|^2+\tfrac14\lambda^2$ by \eqref{eq:poly_post_velocity} and $I'+I_*'=Y$, we obtain from \eqref{eq:E_both}
\begin{equation}
    M(\bv',I')M(\bv_*',I_*')
    = C_M^2\,(I'I_*')^{\frac{\delta}{2}-1}\exp\!\left(-\frac{|\bh-\bu|^2+E}{\theta}\right)
    \leqslant C_M^2\,(I'I_*')^{\frac{\delta}{2}-1}\,e^{-E/\theta},
    \label{eq:MM_decay}
\end{equation}
where $C_M=\rho\,(2\pi\theta)^{-3/2}\theta^{-\delta/2}/\Gamma(\delta/2)$. Hence the contribution of the collisions with $\lambda>S$ decays like $e^{-S^2/(4\theta)}$ as $S$ increases. The same argument applies to the collisions with $g>S$ in both the gain and the loss terms, for which $E>\tfrac14 S^2$ since $U\geqslant0$. Hence $\Lambda=S$ can be taken for polyatomic gases as well when $S$ is sufficiently large.

With $\Lambda=S$, we use the same dealiasing condition as in the monatomic case \cite{pareschi2000numerical,gamba2017fast},
\begin{equation}
    L\geqslant \frac{3+\sqrt2}{4}\,S .
    \label{eq:dealias}
\end{equation}
The resulting $g$-integral on $[0,S]$ is discretized by a quadrature rule $\{w_l^g,g_l\}_{l=1}^{N_g}$, which enters only the offline coefficients.

The integrals with respect to $\lambda$ and $\boldsymbol{\sigma}$ are discretized using numerical quadratures,
\begin{equation}
    \int_0^\Lambda \varphi(\lambda)\,d\lambda \approx \sum_{q=1}^{N_\lambda} w_q^\lambda\,\varphi(\lambda_q),
    \qquad
    \int_{\bbS^2} \varphi(\boldsymbol{\sigma})\,d\boldsymbol{\sigma} \approx \sum_{p=1}^{M} w_p^\sigma\,\varphi(\boldsymbol{\sigma}_p),
    \label{eq:quad_ls}
\end{equation}
with weights $w_q^\lambda,w_p^\sigma$ and nodes $\lambda_q\in[0,\Lambda]$, $\boldsymbol{\sigma}_p\in\bbS^2$.

With these approximations, the Laguerre--Fourier coefficients of the collision term are
\begin{equation}
\begin{aligned}
    \widehat \mQ_{a,\boldsymbol{k}}
    = {}& \frac{1}{(2L)^3}
    \sum_{b,d=0}^{K-1}\sum_{q=1}^{N_\lambda}\sum_{p=1}^{M}
    w_q^\lambda w_p^\sigma\,
    W_{a,b,d}(\lambda_q,\boldsymbol{\sigma}_p,\boldsymbol{\xi}_{\boldsymbol{k}})\,
    \mathcal C_{b,d,\boldsymbol{k}}^{q,p} \\
    &- \frac{1}{(2L)^3}
    \sum_{d=0}^{K-1}(-1)^d\,2^{\alpha+1}
    \sum_{\boldsymbol{m}\in\mathbb K_N}
    \widehat{\mathcal B}_{\boldsymbol{m}-\boldsymbol{k}}\,
    \widehat f_{a,\boldsymbol{m}}\,\widehat f_{d,\boldsymbol{k}-\boldsymbol{m}} ,
\end{aligned}
    \label{eq:Q_disc}
\end{equation}
where $\widehat{\mathcal B}_{\boldsymbol{k}}=\widehat{\mathcal B}(\boldsymbol{\xi}_{\boldsymbol{k}})$ is the loss kernel \eqref{eq:ls_ker} evaluated at the Fourier modes, which depends only on $|\boldsymbol{\xi}_{\boldsymbol{k}}|$, and the multi-indices are understood cyclically. The weights $W_{a,b,d}(\lambda_q,\boldsymbol{\sigma}_p,\boldsymbol{\xi}_{\boldsymbol{k}})$ and $\widehat{\mathcal B}_{\boldsymbol{k}}$ are independent of the distribution function and are precomputed offline, with their explicit quadratures given in App.~\ref{app:offline_coefficients}. The discrete convolution $\mathcal C_{b,d,\boldsymbol{k}}^{q,p}$ appearing in \eqref{eq:Q_disc} is
\begin{equation}
\begin{aligned}
    \mathcal C_{b,d,\boldsymbol{k}}^{q,p}
    =
    \sum_{\boldsymbol{m}\in\mathbb K_N}
    e^{\frac{i}{2}\lambda_q\boldsymbol{\sigma}_p\cdot\boldsymbol{\xi}_{\boldsymbol{m}}}
    \widehat f_{b,\boldsymbol{m}}\,
    e^{-\frac{i}{2}\lambda_q\boldsymbol{\sigma}_p\cdot\boldsymbol{\xi}_{\boldsymbol{k}-\boldsymbol{m}}}
    \widehat f_{d,\boldsymbol{k}-\boldsymbol{m}} .
\end{aligned}
    \label{eq:C_disc}
\end{equation}

For fixed $b,d,q,p$, the convolution $\mathcal C_{b,d,\boldsymbol{k}}^{q,p}$ in \eqref{eq:C_disc} is evaluated for all $\boldsymbol{k}\in\mathbb K_N$ by fast
Fourier transforms in $O(N^3\log N)$ operations. Thus, forming all these convolutions costs
$O(K^2 N_\lambda M N^3\log N)$.

To evaluate the gain term in \eqref{eq:Q_disc}, the convolutions are multiplied by the weights $W_{a,b,d}(\lambda_q,\boldsymbol{\sigma}_p,\boldsymbol{\xi}_{\boldsymbol{k}})$ and summed over $b$, $d$, $q$, and $p$. Evaluated for all $a$ and $\boldsymbol{k}$, these sums require $O(K^3 N_\lambda M N^3)$ operations. The loss term can be evaluated in $O(KN^3\log N)$ operations by first combining the Laguerre modes. The procedure is summarized in
Alg.~\ref{alg:fast_poly}.

Consequently, when $N_\lambda=O(N)$, the overall online computational cost, in which the $O(KN^3\log N)$ cost of the loss term is of lower order, is
\begin{equation}
    C_{\mathrm{fast}}
    = O\!\left(K^2 M N^4\log N + K^3 M N^4\right).
    \label{eq:cost_fast}
\end{equation}
The first term comes from the fast Fourier transforms and the second from the summations with the weights $W_{a,b,d}$, and the dependence on $N$ is the same as in the monatomic fast spectral method \cite{gamba2017fast}. For the VHS-type kernel \eqref{eq:vhs_iso}, $W_{a,b,d}$ is independent of $\boldsymbol{\sigma}$, so that the convolutions can be summed over $p$ before the multiplication by $W_{a,b,d}$, and the online cost reduces to
\begin{equation}
    C_{\mathrm{fast}}^{\mathrm{VHS}}
    = O\!\left(K^2 M N^4\log N + K^3 N^4\right).
    \label{eq:cost_fast_vhs}
\end{equation}

We finally discuss the conservation properties of the discrete collision operator. As in the monatomic case \cite{gamba2017fast,yin2025fast}, the Fourier spectral method does not preserve momentum and energy. Conservation of mass, which by \eqref{eq:Phi_int} reads
\begin{equation}
    \sum_{a=0}^{K-1}(-1)^a\,\widehat \mQ_{a,\boldsymbol 0}=0 .
    \label{eq:mass_disc}
\end{equation}
would hold exactly if the gain and loss terms at $\boldsymbol k=\boldsymbol 0$ were evaluated with the same quadrature. In the scheme \eqref{eq:Q_disc} they are evaluated with different quadratures, so their contributions do not cancel exactly, and \eqref{eq:mass_disc} is satisfied only up to an error of the same order as the discretization error.

\begin{algorithm}[t]
\caption{Fast Laguerre--Fourier evaluation of the collision coefficients $\widehat \mQ_{a,\boldsymbol{k}}$}
\label{alg:fast_poly}
\begin{algorithmic}[1]
\State \textbf{Offline:} choose $S$, $\Lambda$, $L$ and the quadrature rules; precompute the weights $W_{a,b,d}(\lambda_q,\boldsymbol{\sigma}_p,\boldsymbol{\xi}_{\boldsymbol{k}})$ and the loss kernel $\widehat{\mathcal B}_{\boldsymbol{k}}$.
\State \textbf{Input:} Fourier coefficients $\widehat f_{b,\boldsymbol{k}}$, $b=0,\dots,K-1$, $\boldsymbol{k}\in\mathbb K_N$.
\State $\widehat \mQ_{a,\boldsymbol{k}}^+\gets0$ for all $a$ and $\boldsymbol{k}$
\For{$q=1,\dots,N_\lambda$}
  \For{$p=1,\dots,M$}
    \For{$b=0,\dots,K-1$}
      \State $u_b^{\pm}\gets \mathcal F_N^{-1}\big[ e^{\pm\frac{i}{2}\lambda_q\boldsymbol{\sigma}_p\cdot\boldsymbol{\xi}_{\boldsymbol{k}}}\, \widehat f_{b,\boldsymbol{k}}\big]$ \Comment{$2K$ inverse FFTs}
    \EndFor
    \For{$b,d=0,\dots,K-1$}
      \State $\mathcal C_{b,d,\boldsymbol{k}}^{q,p}\gets (2L)^3\,\mathcal F_N\big[u_b^{+}\,u_d^{-}\big]_{\boldsymbol{k}}$ \Comment{$K^2$ forward FFTs}
      \For{$a=0,\dots,K-1$}
        \State $\widehat \mQ_{a,\boldsymbol{k}}^+\gets \widehat \mQ_{a,\boldsymbol{k}}^+ + \dfrac{w_q^\lambda w_p^{\sigma}}{(2L)^3}\,W_{a,b,d}(\lambda_q,\boldsymbol{\sigma}_p,\boldsymbol{\xi}_{\boldsymbol{k}})\, \mathcal C_{b,d,\boldsymbol{k}}^{q,p}$
      \EndFor
    \EndFor
  \EndFor
\EndFor
\State $\widehat n_{\boldsymbol{k}}\gets 2^{\alpha+1}\sum_{d=0}^{K-1} (-1)^d\,\widehat f_{d,\boldsymbol{k}}$ \Comment{$O(KN^3)$}
\For{$a=0,\dots,K-1$}
  \State $\widehat \mQ_{a,\boldsymbol{k}}^-\gets \mathcal F_N\big[\mathcal F_N^{-1}[\widehat f_{a,\boldsymbol{k}}]\;\mathcal F_N^{-1}[\widehat{\mathcal B}_{\boldsymbol{k}}\,\widehat n_{\boldsymbol{k}}]\big]_{\boldsymbol{k}}$ \Comment{$K$ convolutions}
\EndFor
\State \textbf{Output:} $\widehat \mQ_{a,\boldsymbol{k}}\gets \widehat \mQ_{a,\boldsymbol{k}}^+-\widehat \mQ_{a,\boldsymbol{k}}^-$, $a=0,\dots,K-1$, $\boldsymbol{k}\in\mathbb K_N$.
\end{algorithmic}
\end{algorithm}

\section{Numerical results}
\label{sec:numerical}

In this section, we use numerical experiments to evaluate the accuracy and robustness of the proposed fast Laguerre--Fourier spectral method for the polyatomic Boltzmann equation \eqref{eq:bolt}. We test the spatially homogeneous problem in Sec.~\ref{subsec:homogeneous} and the one-dimensional steady problem in Sec.~\ref{subsec:flows}. In all the tests we take $\delta=2$ internal degrees of freedom and $B_0=1/(4\pi)$ in the kernel \eqref{eq:vhs_iso}; the exponent $\zeta$ is specified for each test.

\subsection{Spatially homogeneous Boltzmann equation}
\label{subsec:homogeneous}

The tests in this subsection solve the spatially homogeneous equation
\begin{equation}
    \partial_t f = \mQ(f,f), \qquad f(0,\bv,I)=f_0(\bv,I),
    \label{eq:homogeneous}
\end{equation}
where $f_0$ is the initial condition; this equation retains only the collision term of \eqref{eq:bolt}. Unless otherwise stated, the collision discretization uses the relative-speed cutoff $S=14$, $L=(3+\sqrt2)S/4$, $\Lambda=S$, $N_\lambda=N$, and a Lebedev rule with $M=26$ angular nodes, and the time derivative is discretized using the classical fourth-order Runge--Kutta scheme with time step $\Delta t=0.1$. The offline coefficients in App.~\ref{app:offline_coefficients} are computed with $N_g=N$ quadrature points in $g$, $N_r=N_\tau=16$ points for $G_{b,d}$ and $\mathcal P_a$, and an $N_Y=48$ point Gauss--Laguerre rule in $Y$. The scaling parameter of the Laguerre functions \eqref{eq:psi_b} is fixed at $\beta=1$ in all the tests.

\subsubsection{Isotropic initial condition}
\label{subsec:accuracy}

We first consider the isotropic initial condition of BKW type,
\begin{equation}
\begin{aligned}
    f_0(\bv,I)
    ={}&
    \frac{I^{\frac{\delta}{2}-1}}
    {(2\pi\kappa)^{3/2}\kappa^{\delta/2}\Gamma(\delta/2)}
    \exp\!\left(-\frac{\frac12|\bv|^2+I}{\kappa}\right) \\
    &\times
    \left(
    \frac{(D+2)\kappa-D}{2\kappa}
    +\frac{1-\kappa}{2\kappa^2}\,
    \big(|\bv|^2+2I\big)
    \right),
\end{aligned}
    \label{eq:poly_bkw}
\end{equation}
where $D=3+\delta$ and $\kappa=0.75$; in this test the kernel is of Maxwell type, $\zeta=0$. The polynomial factor makes $f_0$ non-Maxwellian in the velocity and in the internal energy at the same time, so that the test exercises the two discretizations together.

\begin{figure}[!hptb]
\centering
\subfigure[$L^2$ error for different $N$ with $K=20$ fixed.]{%
    \includegraphics[width=0.48\textwidth]{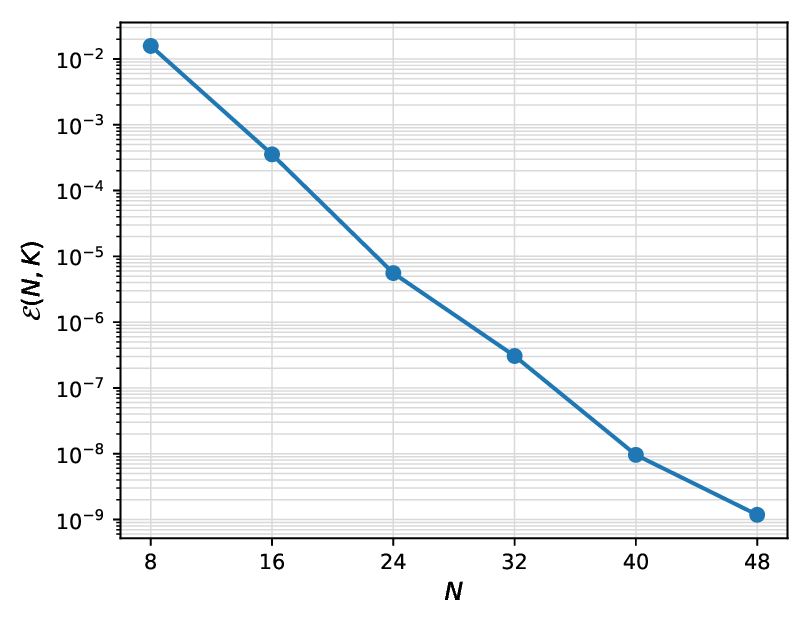}%
    \label{fig:conv_N}%
}
\subfigure[$L^2$ error for different $K$ with $N=96$ fixed.]{%
    \includegraphics[width=0.48\textwidth]{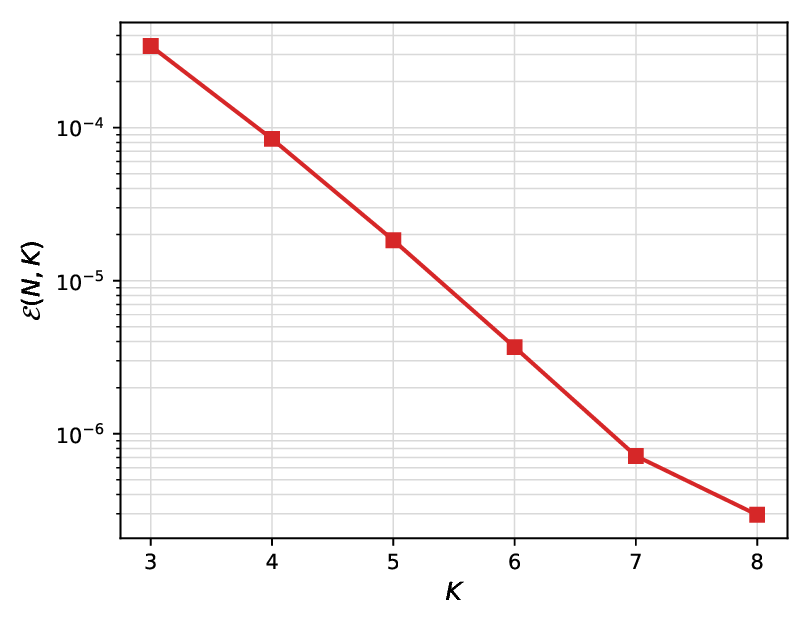}%
    \label{fig:conv_K}%
}
\caption{(Isotropic initial condition in Sec.~\ref{subsec:accuracy}) Spectral convergence of the collision operator $\mQ(f_0,f_0)$ in $N$ and $K$: the $L^2$ error $\mathcal E(N,K)$ on a semi-log scale.}
\label{fig:conv}
\end{figure}

We first measure the accuracy of the collision discretization itself. The collision term $\mQ(f_0,f_0)$ is evaluated on a sequence of grids and compared with a high-resolution reference evaluation through the $L^2$ error
\begin{equation}
    \mathcal E(N,K)
    =\big\|\,\mQ_{N,K}(f_0,f_0)-\mQ_{\mathrm{fine}}(f_0,f_0)\,\big\|_2 ,
    \label{eq:conv_error}
\end{equation}
where the reference $\mQ_{\mathrm{fine}}$ is computed with $(N,K)=(96,20)$. The two directions are treated separately: in Fig.~\ref{fig:conv_N} the number of Laguerre modes is fixed at $K=20$ and $N$ is increased from $8$ to $48$, and in Fig.~\ref{fig:conv_K} the number of Fourier modes is fixed at $N=96$ and $K$ is increased from $3$ to $8$. In both cases the error falls on a straight line in the semi-log plot, which indicates the spectral convergence of the Fourier and Laguerre discretizations.

We next measure the cost of the method. Tab.~\ref{tab:cost} lists the CPU time of a single evaluation of $\mQ(f_0,f_0)$, for $N$ from $8$ to $96$ at $K=20$ fixed, and for $K$ from $3$ to $10$ at $N=96$ fixed. Following Alg.~\ref{alg:fast_poly}, the time is split into the offline precomputation of the coefficients and the online evaluation. The offline cost grows much more slowly than the online one and is paid once for an entire simulation.

The online time is used to verify the estimate \eqref{eq:cost_fast_vhs}. The order reported in Tab.~\ref{tab:cost} is that of the online time between two consecutive resolutions in Tab.~\ref{tab:cost},
\begin{equation}
    p_i = \frac{\log(t_i/t_{i-1})}{\log(x_i/x_{i-1})}, \qquad x=N \ \text{or} \ K ,
    \label{eq:measured_order}
\end{equation}
so that $p_i$ is the exponent of a local power law $t\propto x^{p_i}$.

When $N$ is increased at $K=20$ fixed, the measured order grows from about $1$ on the coarsest grids to more than $4$ on the finest. At small $N$ the evaluation is dominated by costs that do not depend on the grid size, such as the overhead of each transform call, so that the measured time reflects little more than the number of transforms. As $N$ grows these fixed costs become negligible against the work of the transforms themselves, and the order reaches $4.12$ and $4.07$ at the last two steps, which lie between the order $4$ of the second term in \eqref{eq:cost_fast_vhs} and the local order $4+1/\ln N\approx4.2$ of its first term $K^2MN^4\log N$.

Increasing $K$ at $N=96$ fixed shows no such transition, since increasing $K$ changes the number of transforms but not their size, so that the asymptotic regime is reached at once. The measured order fluctuates between $2.08$ and $2.81$, that is, between the $K^2$ of the convolutions and the $K^3$ of the contraction in \eqref{eq:cost_fast_vhs}.

\begin{table}[!htbp]
\centering
\caption{(Isotropic initial condition in Sec.~\ref{subsec:accuracy}) CPU time of a single evaluation of the collision operator $\mQ(f_0,f_0)$, split into the offline computation of the coefficients and the online evaluation. The order is that of the online time, defined in \eqref{eq:measured_order}.}
\label{tab:cost}
\setlength{\tabcolsep}{8pt}
\begin{tabular}{cccc@{\qquad\qquad}cccc}
\toprule
\multirow{2}{*}{$N$} & \multicolumn{3}{c}{CPU time (s) for $K=20$} &
\multirow{2}{*}{$K$} & \multicolumn{3}{c}{CPU time (s) for $N=96$}\\
\cmidrule(lr){2-4}\cmidrule(lr){6-8}
 & offline & online & order & & offline & online & order\\
\midrule
  8 &   61.7 &    10.5 &    -- &  3 &   28.0 &   26.7 &    --\\
 16 &  115.9 &    19.8 &  0.92 &  4 &   38.9 &   53.8 &  2.44\\
 24 &  170.7 &    29.0 &  0.94 &  5 &   54.1 &   89.2 &  2.27\\
 32 &  224.6 &    38.7 &  1.00 &  6 &   72.2 &  134.2 &  2.24\\
 40 &  280.2 &    69.4 &  2.62 &  7 &   94.8 &  193.5 &  2.37\\
 48 &  332.8 &   152.9 &  4.33 &  8 &  119.9 &  255.4 &  2.08\\
 64 &  449.2 &   500.5 &  4.12 &  9 &  144.5 &  355.4 &  2.81\\
 96 &  668.5 &  2609.6 &  4.07 & 10 &  176.2 &  461.7 &  2.48\\
\bottomrule
\end{tabular}
\end{table}

\begin{figure}[!hptb]
\centering
\subfigure[internal-energy cut $f(\bv_0,I)$]{\includegraphics[width=0.46\textwidth]{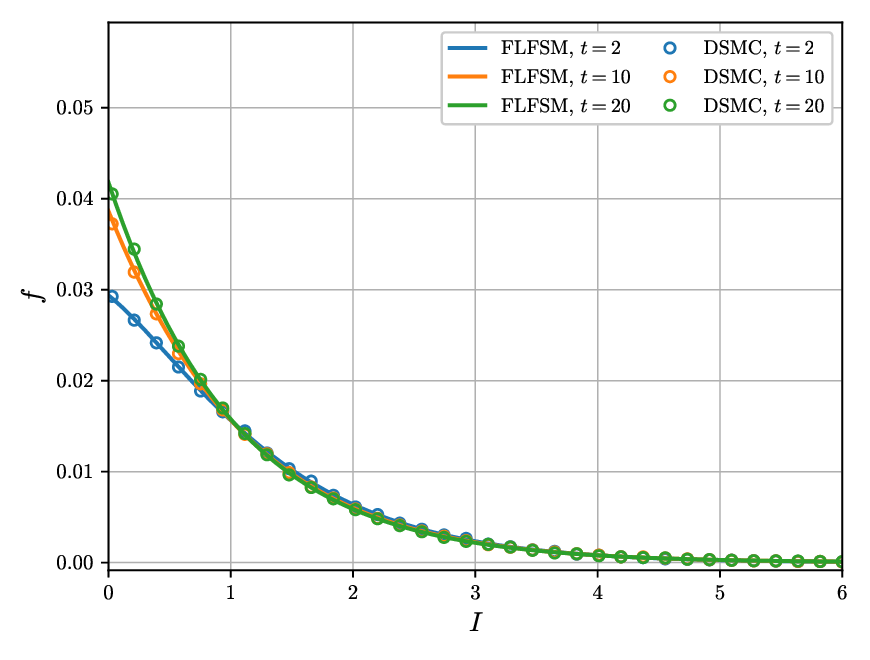}}
\hfill
\subfigure[velocity cut $f(v_1,0,0,I_0)$ at $I_0=0.5$]{\includegraphics[width=0.46\textwidth]{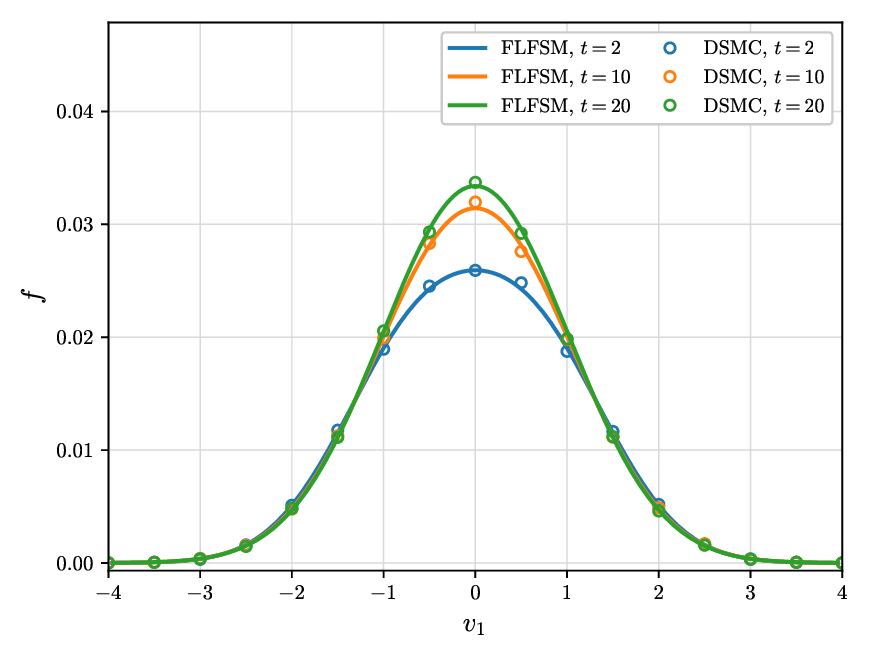}}
\caption{(Isotropic initial condition in Sec.~\ref{subsec:accuracy}) Relaxation from the isotropic initial condition \eqref{eq:poly_bkw} at $t=2,10,20$: (a) the internal-energy cut at $\bv_0=(0.368,0,0)$ and (b) the velocity cut $f(v_1,0,0,I_0)$ at $I_0=0.5$. The solid lines are the FLFSM solutions and the circles are the DSMC solutions.}
\label{fig:bkw_lines}
\end{figure}

We then solve the homogeneous equation \eqref{eq:homogeneous} from the same initial condition with $N=72$ and $K=10$ modes, up to $t=20$. Since the solution is four-dimensional, we show only two one-dimensional cuts, one along the internal energy and one along a velocity component. Fig.~\ref{fig:bkw_lines}(a) shows the evolution of the internal-energy cut $f(\bv_0,I)$ at $\bv_0=(0.368,0,0)$ and Fig.~\ref{fig:bkw_lines}(b) that of the velocity cut $f(v_1,0,0,I_0)$ at $I_0=0.5$, both at $t=2,10,20$ and compared with a DSMC solution. The two solutions essentially coincide at all three times, which supports the accuracy of the FLFSM, and the curves at $t=10$ and $t=20$ lie close to each other, showing that the distribution function is approaching equilibrium.

\subsubsection{Anisotropic initial condition}
\label{subsec:homogeneous_relaxation}
To further test the performance of the FLFSM, we consider two anisotropic initial conditions adapted from \cite{yin2025fast}, to which an internal-energy dependence has been added. The first one is a smooth superposition of four Maxwellians
\begin{equation}
\begin{aligned}
    f^{(1)}(\bv,I)=\frac{1}{4\,(2\pi\theta_0)^{3/2}\theta_0}\bigg[
    &\exp\!\Big(-\frac{(v_1+u_0)^2+v_2^2+v_3^2}{2\theta_0}-\frac{I}{\theta_0}\Big)\\
    &+\exp\!\Big(-\frac{(v_1-u_0)^2+v_2^2+v_3^2}{2\theta_0}-\frac{I}{\theta_0}\Big)\\
    &+\exp\!\Big(-\frac{v_1^2+(v_2+u_0)^2+v_3^2}{2\theta_0}-\frac{I}{\theta_0}\Big)\\
    &+\exp\!\Big(-\frac{v_1^2+(v_2-u_0)^2+v_3^2}{2\theta_0}-\frac{I}{\theta_0}\Big)\bigg],
\end{aligned}
\label{eq:homog_ic1}
\end{equation}
where $u_0=\sqrt2$ and $\theta_0=\tfrac35$. Here and in \eqref{eq:homog_ic2} and \eqref{eq:wall_bc} below, we have used $\delta=2$, so that $I^{\delta/2-1}/\Gamma(\delta/2)=1$.

The second initial condition is discontinuous across $v_1=0$, which makes it a demanding test for a spectral discretization in velocity,
\begin{equation}
    f^{(2)}(\bv,I)=\begin{cases}\dfrac{2^{3/4}(2-\sqrt2)}{\pi^{3/2}}\,\exp\!\Big(-\frac{|\bv|^2}{\sqrt2}-\sqrt2\,I\Big), & v_1>0,\\[8pt]\dfrac{2^{3/4}(2-\sqrt2)}{8\pi^{3/2}}\,\exp\!\Big(-\frac{|\bv|^2}{2\sqrt2}-\frac{I}{\sqrt2}\Big), & v_1<0.\end{cases}
    \label{eq:homog_ic2}
\end{equation}
By construction, both $f^{(1)}$ and $f^{(2)}$ have density $\rho=1$, macroscopic velocity $\bu=\bz$, and total temperature $\theta=1$. The kernel is again of Maxwell type, $\zeta=0$.

\begin{figure}[!hptb]
\centering
\subfigure[$t=0.5$]{\includegraphics[width=0.32\textwidth]{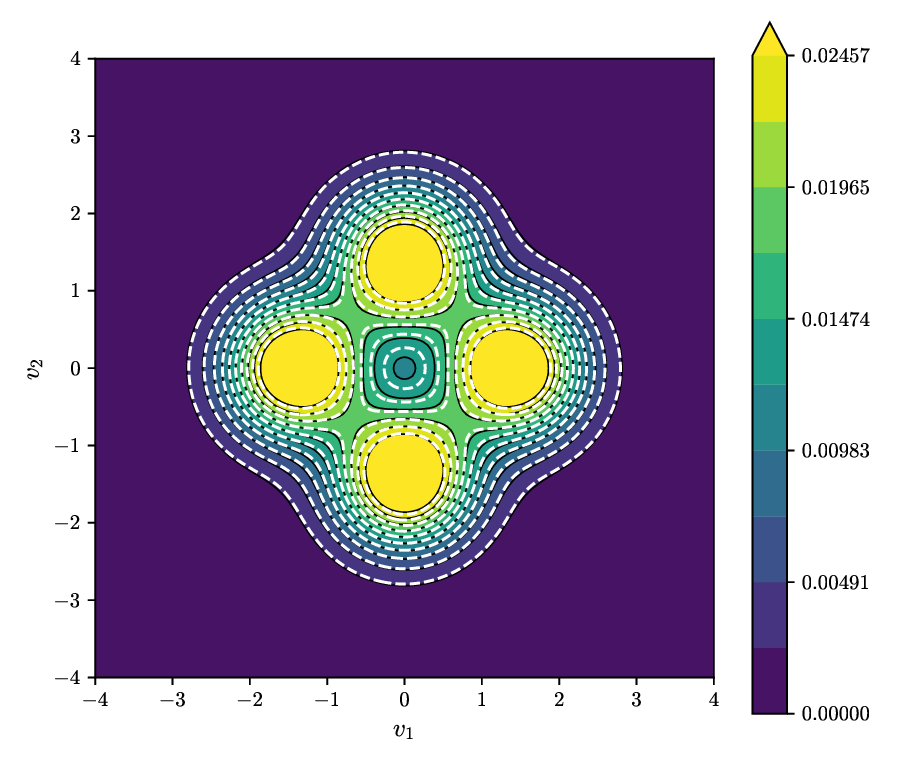}}
\subfigure[$t=4$]{\includegraphics[width=0.32\textwidth]{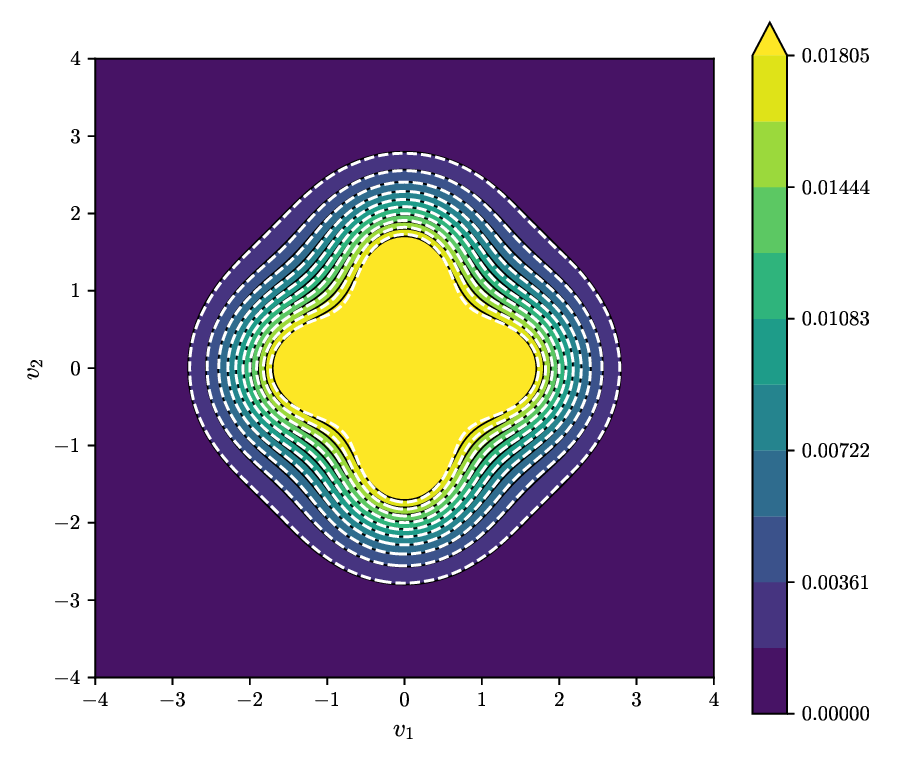}}
\subfigure[$t=12$]{\includegraphics[width=0.32\textwidth]{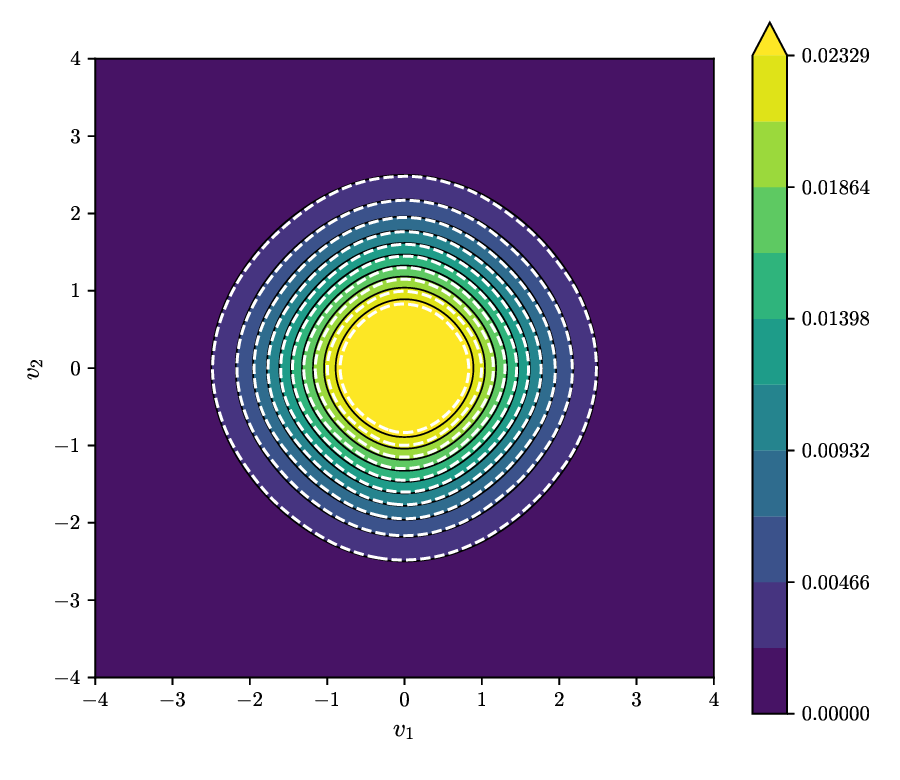}}
\caption{(Anisotropic initial conditions in Sec.~\ref{subsec:homogeneous_relaxation}) Comparison of the velocity-space slices at $v_3=0$, $I=0.5$ between the FLFSM and DSMC for the initial condition \eqref{eq:homog_ic1} at $t=0.5,4,12$. The filled contours with black lines represent the FLFSM solutions, and the gray dashed lines denote the DSMC solutions.}
\label{fig:homog_ic1}
\end{figure}

\begin{figure}[!hptb]
\centering
\subfigure[$t=0.5$]{\includegraphics[width=0.32\textwidth]{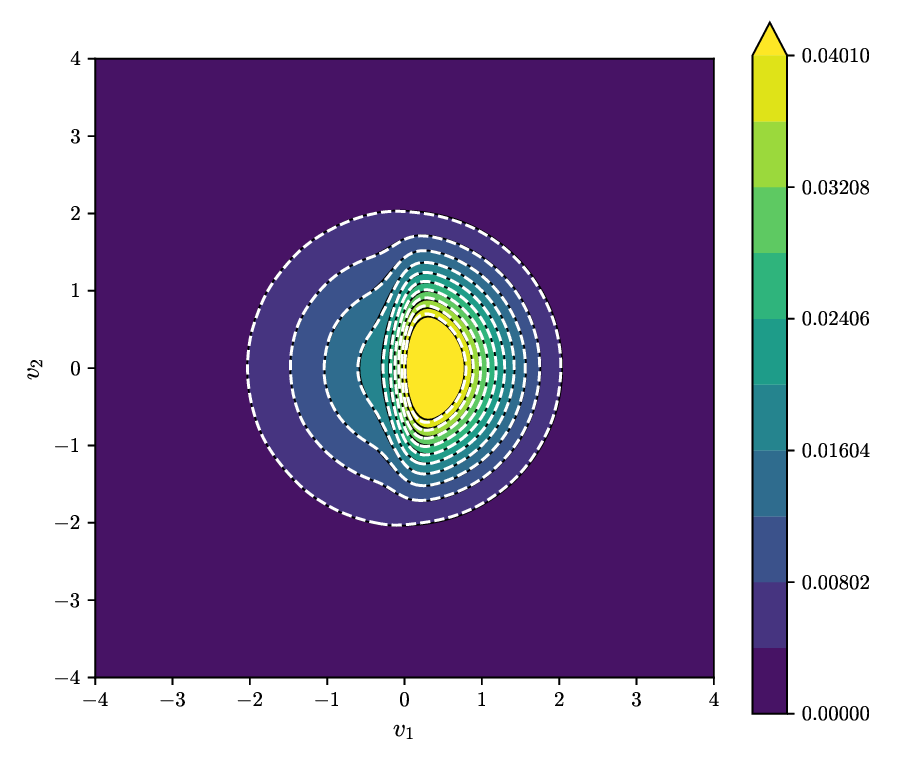}}
\subfigure[$t=4$]{\includegraphics[width=0.32\textwidth]{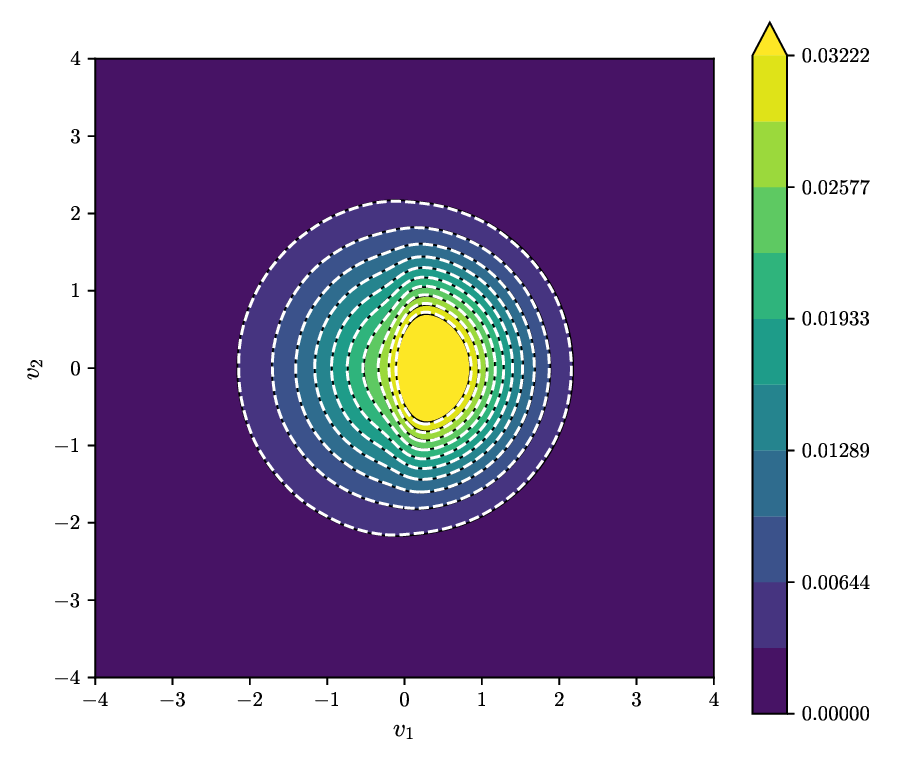}}
\subfigure[$t=12$]{\includegraphics[width=0.32\textwidth]{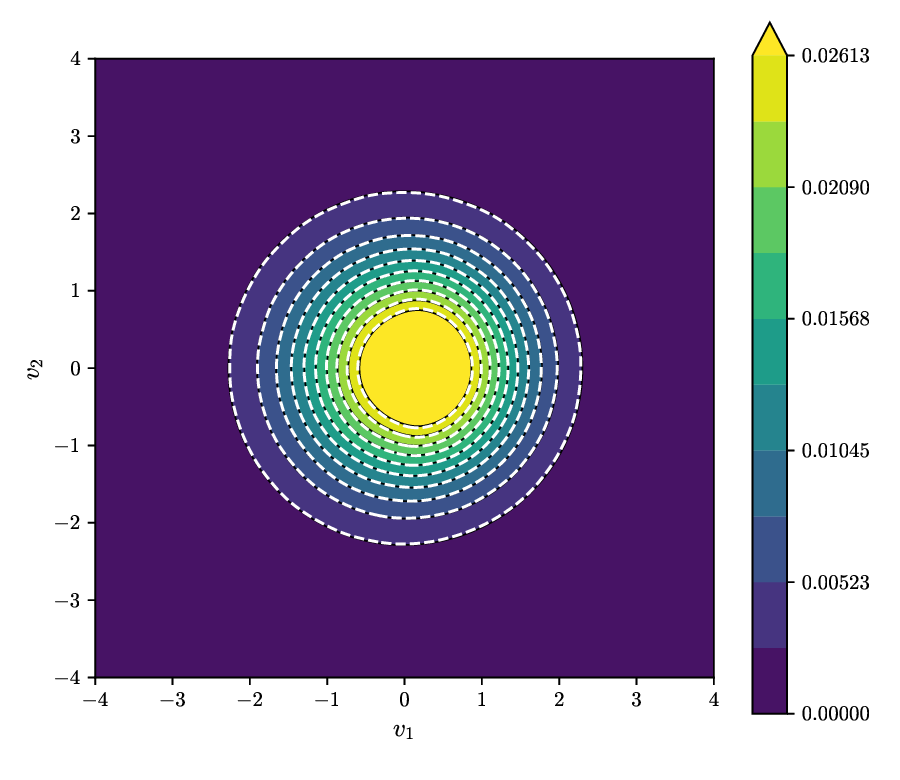}}
\caption{(Anisotropic initial conditions in Sec.~\ref{subsec:homogeneous_relaxation}) Comparison of the velocity-space slices at $v_3=0$, $I=0.5$ between the FLFSM and DSMC for the initial condition \eqref{eq:homog_ic2} at $t=0.5,4,12$. The filled contours with black lines represent the FLFSM solutions, and the gray dashed lines denote the DSMC solutions.}
\label{fig:homog_ic2}
\end{figure}

Then, we solve the homogeneous equation \eqref{eq:homogeneous} with the initial conditions \eqref{eq:homog_ic1} and \eqref{eq:homog_ic2}. With $N=72$ and $K=10$ modes, the solution is evolved up to $t=12$. Since the solution is four-dimensional, we plot only the velocity-space slice at $v_3=0$ and $I=0.5$. Figs.~\ref{fig:homog_ic1} and \ref{fig:homog_ic2} show the FLFSM solution for $f^{(1)}$ and $f^{(2)}$, respectively, at $t=0.5,4,12$, together with the DSMC solution. For both initial conditions the FLFSM and DSMC solutions lie on top of each other at all three times, including where $f^{(2)}$ is discontinuous and a truncated Fourier series cannot represent it exactly. During the evolution, the four peaks of $f^{(1)}$ and the jump of $f^{(2)}$ are smoothed out, and both solutions become isotropic and approach the equilibrium Maxwellian.

\subsubsection{Two-temperature relaxation test}
\label{subsec:relaxation_dsmc}

A key difference between polyatomic and monatomic gases lies in the exchange between translational and internal energy, which can be observed through the evolution of the translational and internal temperatures. Although no exact solution of the polyatomic Boltzmann equation is known at the level of the distribution function, for the Maxwell kernel $\zeta=0$ the two temperatures obey the closed analytic law derived in App.~\ref{subsec:two_temperature_law}. This allows the FLFSM to be validated at the macroscopic level against an exact solution, and we consider the two-temperature initial condition
\begin{equation}
    f_0(\bv,I) = \frac{I^{\delta/2-1}}{(2\pi\theta_{\mathrm{tr}}^0)^{3/2}(\theta_{\mathrm{int}}^0)^{\delta/2}\Gamma(\delta/2)}
    \exp\!\left(-\frac{|\bv|^2}{2\theta_{\mathrm{tr}}^0}-\frac{I}{\theta_{\mathrm{int}}^0}\right),
    \label{eq:two_temp_initial}
\end{equation}
in which the translational temperature $\theta_{\mathrm{tr}}^0$ and the internal temperature $\theta_{\mathrm{int}}^0$ are different. Starting from \eqref{eq:two_temp_initial}, the two temperatures relax exponentially to a common value,
\begin{equation}
    \theta_{\mathrm{tr}}(t) = \theta_\infty + \big(\theta_{\mathrm{tr}}^0-\theta_\infty\big)e^{-\mu t},
    \qquad
    \theta_{\mathrm{int}}(t) = \theta_\infty + \big(\theta_{\mathrm{int}}^0-\theta_\infty\big)e^{-\mu t},
    \label{eq:twoT_main}
\end{equation}
where the equilibrium temperature $\theta_\infty=(3\theta_{\mathrm{tr}}^0+\delta\theta_{\mathrm{int}}^0)/(3+\delta)$ is fixed by the conservation of energy and the relaxation rate is $\mu=(3+\delta)\nu_0/(2\delta+3)$, with $\nu_0$ given in \eqref{eq:nu0}. For $B_0=1/(4\pi)$ and $\delta=2$ this gives $\mu=4/21$.

\begin{figure}[!hptb]
\centering
\subfigure[Temperature error for different $N$ with $K=20$ fixed.]{%
    \includegraphics[width=0.42\textwidth]{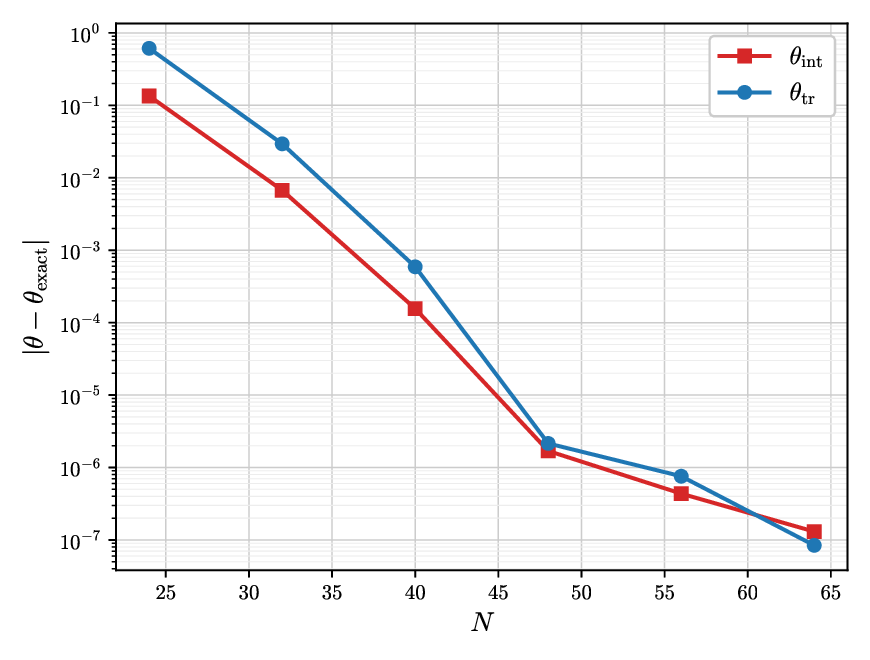}%
    \label{fig:relax_conv_N}%
}
\hfill
\subfigure[Temperature error for different $K$ with $N=80$ fixed.]{%
    \includegraphics[width=0.42\textwidth]{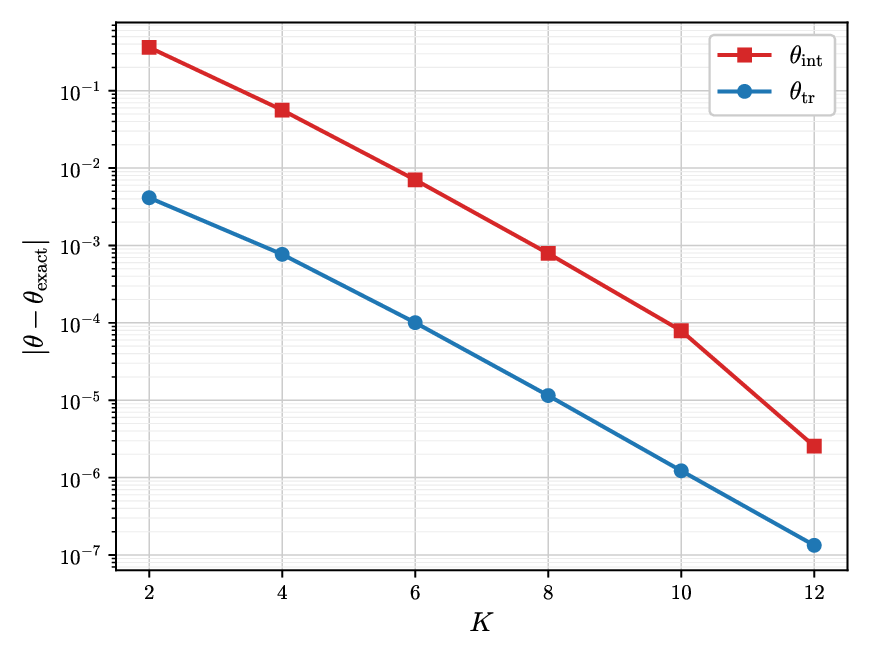}%
    \label{fig:relax_conv_K}%
}
\caption{(Two-temperature relaxation with $\zeta=0$ in Sec.~\ref{subsec:relaxation_dsmc}) Spectral convergence of the macroscopic temperatures at $t=0.5$ in $N$ and $K$: the error with respect to the analytic solution \eqref{eq:twoT_main} on a semi-log scale.}
\label{fig:relax_conv}
\end{figure}

In the following tests, we take $\zeta=0$ and set $\theta_{\mathrm{tr}}^0=1.3$ and $\theta_{\mathrm{int}}^0=0.55$, so that $\theta_\infty=1$. We test the convergence by evolving the homogeneous equation \eqref{eq:homogeneous} with the initial condition \eqref{eq:two_temp_initial} up to $t=0.5$ with $\Delta t=0.025$. The absolute error between the numerical and exact solutions, $|\theta-\theta_{\text{exact}}|$, for the two temperatures is shown in Fig.~\ref{fig:relax_conv}. In the velocity direction, we fix $K=20$ and increase $N$ from $24$ to $64$, with the results shown in Fig.~\ref{fig:relax_conv_N}. In the internal-energy direction, we fix $N=80$ and increase $K$ from $2$ to $12$, with the results shown in Fig.~\ref{fig:relax_conv_K}. In both cases the computed temperatures converge rapidly to the analytic solution \eqref{eq:twoT_main}, with the error decaying linearly on the semi-log scale, confirming the spectral accuracy of the FLFSM at the macroscopic level.

\begin{figure}[!hptb]
\centering
\subfigure[Temperature evolution.]{%
    \includegraphics[width=0.48\textwidth]{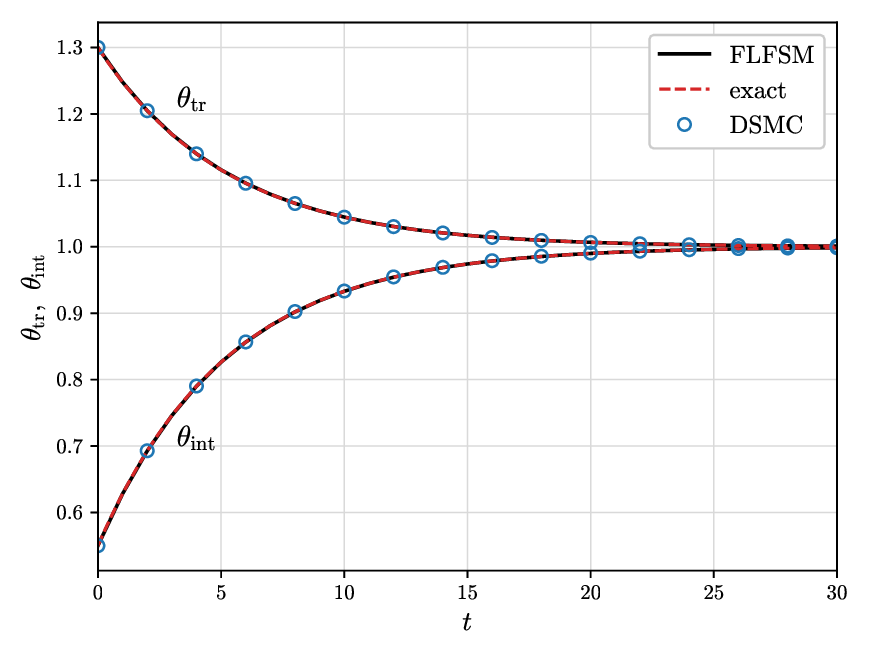}%
    \label{fig:relax_traj}%
}
\subfigure[Decay rate of the temperature difference.]{%
    \includegraphics[width=0.48\textwidth]{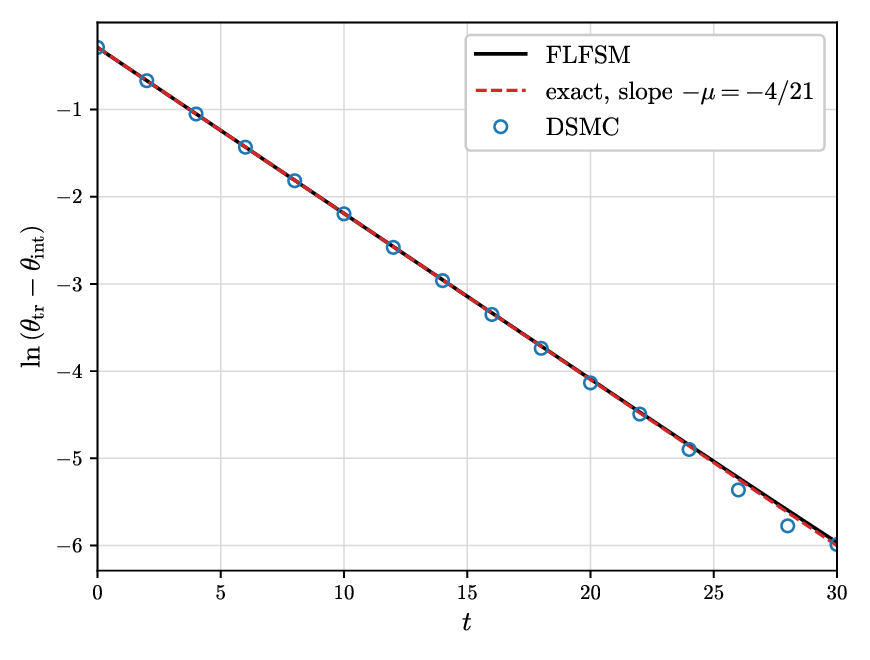}%
    \label{fig:relax_rate}%
}
\caption{(Two-temperature relaxation with $\zeta=0$ in Sec.~\ref{subsec:relaxation_dsmc}) Comparison of the macroscopic temperatures between the FLFSM, DSMC, and the analytic solution \eqref{eq:twoT_main} for the two-temperature initial condition \eqref{eq:two_temp_initial}: (a) the translational and internal temperatures versus time, and (b) the decay of the temperature difference $\theta_{\mathrm{tr}}-\theta_{\mathrm{int}}$ on a semi-log scale, with the analytic slope $-\mu=-4/21$.}
\label{fig:relaxation}
\end{figure}

We then perform a long-time relaxation test, evolving the homogeneous equation \eqref{eq:homogeneous} with the same initial condition \eqref{eq:two_temp_initial} up to $t=30$ with $\Delta t=0.1$ using $N=48$ and $K=10$ modes. Fig.~\ref{fig:relax_traj} shows the translational and internal temperatures against the DSMC results and the analytic solution \eqref{eq:twoT_main}, and Fig.~\ref{fig:relax_rate} shows the decay of the temperature difference $\theta_{\mathrm{tr}}-\theta_{\mathrm{int}}$. The FLFSM, DSMC, and analytic results remain in close agreement over the whole relaxation, with the temperature difference decaying as a single exponential at the analytic rate $\mu=4/21$.

\begin{figure}[!hptb]
\centering
\subfigure[Temperature evolution.]{%
    \includegraphics[width=0.48\textwidth]{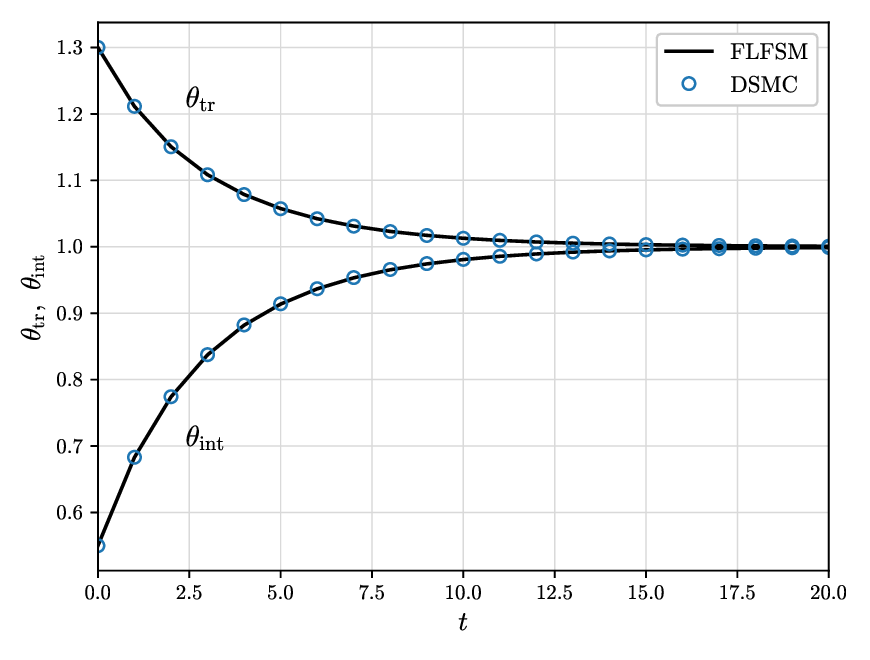}%
    \label{fig:relax_zeta1_traj}%
}
\subfigure[Decay rate of the temperature difference.]{%
    \includegraphics[width=0.48\textwidth]{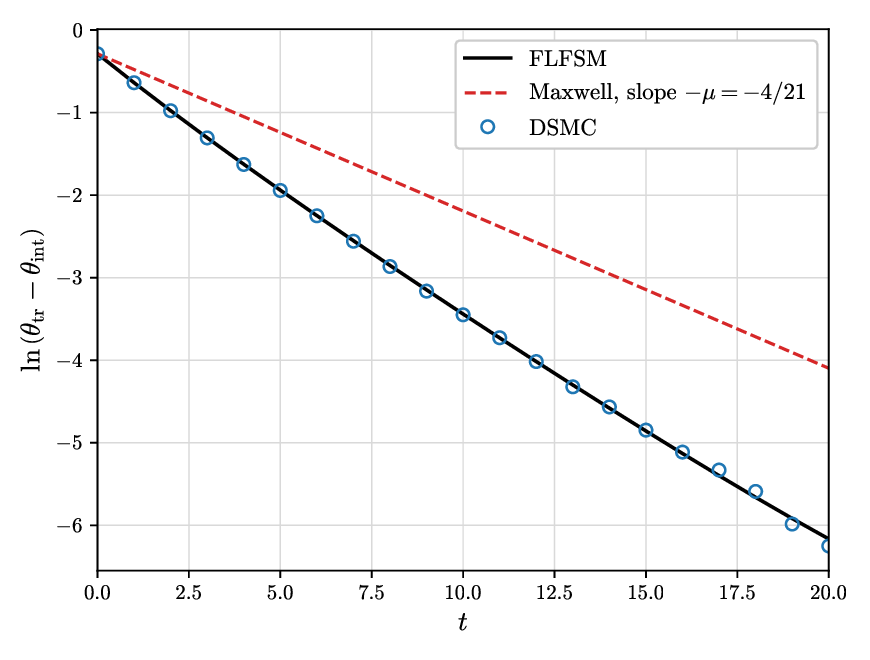}%
    \label{fig:relax_zeta1_rate}%
}
\caption{(Two-temperature relaxation with $\zeta=1$ in Sec.~\ref{subsec:relaxation_dsmc}) Comparison of the macroscopic temperatures between the FLFSM and DSMC for the two-temperature initial condition \eqref{eq:two_temp_initial}: (a) the translational and internal temperatures versus time, and (b) the decay of the temperature difference $\theta_{\mathrm{tr}}-\theta_{\mathrm{int}}$ on a semi-log scale. Since no analytic reference is available for $\zeta=1$, the dashed line in (b) shows the Maxwell law \eqref{eq:twoT_main} for comparison only.}
\label{fig:relax_zeta1}
\end{figure}

To test the method on a more general collision kernel, we repeat the relaxation with the hard-sphere kernel $\zeta=1$, for which the closure of App.~\ref{subsec:two_temperature_law} no longer applies and no analytic law for the two temperatures is available. All other parameters are unchanged, except that the solution is evolved up to $t=20$, and the result is validated against DSMC alone. Fig.~\ref{fig:relax_zeta1_traj} shows the evolution of the two temperatures, and Fig.~\ref{fig:relax_zeta1_rate} the decay of their difference. The FLFSM and DSMC results follow each other throughout, and the gas has essentially reached equilibrium by $t\approx20$. Fig.~\ref{fig:relax_zeta1_rate} also shows the Maxwell law \eqref{eq:twoT_main}, which no longer applies here and is drawn only for comparison. The hard-sphere relaxation is clearly faster, and its decay is no longer a single exponential. The FLFSM therefore applies equally to non-Maxwell kernels such as the hard-sphere kernel, and not only to the exactly solvable Maxwell case.

\subsection{Spatially one-dimensional steady-state problem}
\label{subsec:flows}

In this subsection, numerical experiments are carried out to test the performance of the proposed method on one-dimensional microflows, including the Fourier and Couette flows. Both use the Maxwell kernel, $\zeta=0$. 

\subsubsection{Problem settings}
\label{subsec:flow_setup}

We consider the 1D steady-state Boltzmann equation
\begin{equation}
    v_1\,\partial_x f = \frac{1}{\Kn}\mQ(f,f), \qquad x\in[0,1],
    \label{eq:steady_1d}
\end{equation}
where $v_1$ is the velocity component in the $x$-direction, and $\Kn$ is the Knudsen number measuring the flow rarefaction. The gas is confined between two infinite parallel plates at $x=0$ and $x=1$, and Maxwell's diffuse-reflection boundary conditions \cite{maxwell1878iii} are imposed at both walls. Assuming the solid wall has velocity $\bu^W$ and temperature $\theta^W$, the molecules leaving the wall follow the wall Maxwellian
\begin{equation}
    f(\bv,I) = \frac{\rho^W}{(2\pi \theta^W)^{3/2}\,\theta^W}\exp\!\left(-\frac{\tfrac12|\bv-\bu^W|^2+I}{\theta^W}\right), \qquad \bv\cdot\bn>0,
    \label{eq:wall_bc}
\end{equation}
where $\bn$ is the unit normal pointing into the gas and the wall density $\rho^W$ is fixed by the zero net mass flux at the wall,
\begin{equation}
    \rho^W=\sqrt{\frac{2\pi}{\theta^W}}\int_{\bv\cdot\bn<0}|\bv\cdot\bn|\,f\,d\bv\,dI .
    \label{eq:wall_rho}
\end{equation}

We solve the steady equation \eqref{eq:steady_1d} by source iteration. Let $f^{(\ell)}$ be the result after the $\ell$-th iteration. Then $f^{(\ell+1)}$ is obtained from
\begin{equation}
    v_1\,\partial_x f^{(\ell+1)} + \frac{\nu}{\Kn}\,f^{(\ell+1)}
    = \frac{1}{\Kn}\Big(\mQ\big(f^{(\ell)},f^{(\ell)}\big) + \nu\,f^{(\ell)}\Big),
    \qquad \ell = 0,1,\dots .
    \label{eq:source_iteration}
\end{equation}
Here $\nu=\nu_0=4/15$ is taken as the collision frequency \eqref{eq:nu0} of Maxwell molecules with $\delta=2$. 

We divide $[0,1]$ into $N_x$ uniform cells of size $\Delta x=1/N_x$, with centers $x_i=(i-\tfrac12)\Delta x$ for $i=1,\dots,N_x$, and discretize the velocity and the internal energy by the FLFSM of Sec.~\ref{sec:fast_fourier}. Let $f^{(\ell)}_{a,i,\boldsymbol k}$ denote the $a$-th internal-energy coefficient at the spatial node $x_i$ and the velocity node $\bv_{\boldsymbol k}$, and let $\mQ^{(\ell)}_{a,i,\boldsymbol k}$ be the collision coefficient evaluated from $\{f^{(\ell)}_{a,i,\boldsymbol k}\}$ by Alg.~\ref{alg:fast_poly}. Since the transport term acts independently on each Laguerre mode and each velocity node, the first-order upwind discretization of \eqref{eq:source_iteration} reads
\begin{equation}
\begin{aligned}
    v_{\boldsymbol k,1}\,\frac{f^{(\ell+1)}_{a,i,\boldsymbol k}-f^{(\ell+1)}_{a,i-1,\boldsymbol k}}{\Delta x}
    + \frac{\nu}{\Kn}\,f^{(\ell+1)}_{a,i,\boldsymbol k}
    &= \frac{1}{\Kn}\Big(\mQ^{(\ell)}_{a,i,\boldsymbol k} + \nu\,f^{(\ell)}_{a,i,\boldsymbol k}\Big),
    && v_{\boldsymbol k,1}\geqslant 0,\\
    v_{\boldsymbol k,1}\,\frac{f^{(\ell+1)}_{a,i+1,\boldsymbol k}-f^{(\ell+1)}_{a,i,\boldsymbol k}}{\Delta x}
    + \frac{\nu}{\Kn}\,f^{(\ell+1)}_{a,i,\boldsymbol k}
    &= \frac{1}{\Kn}\Big(\mQ^{(\ell)}_{a,i,\boldsymbol k} + \nu\,f^{(\ell)}_{a,i,\boldsymbol k}\Big),
    && v_{\boldsymbol k,1}< 0 .
\end{aligned}
    \label{eq:upwind}
\end{equation}
For a velocity with $v_{\boldsymbol k,1}\geqslant 0$, \eqref{eq:upwind} determines $f^{(\ell+1)}_{a,i,\boldsymbol k}$ explicitly from $f^{(\ell+1)}_{a,i-1,\boldsymbol k}$, so the new iterate is obtained by sweeping from $i=1$ to $i=N_x$; for $v_{\boldsymbol k,1}<0$ the sweep runs from $i=N_x$ to $i=1$. The wall conditions \eqref{eq:wall_bc}--\eqref{eq:wall_rho} are imposed through the ghost cells $i=0$ and $i=N_x+1$: the wall Maxwellian \eqref{eq:wall_bc} is expanded in the same Laguerre basis \eqref{eq:f_exp}, and the wall density \eqref{eq:wall_rho} is computed from the part of $f^{(\ell)}$ in the adjacent cell that is incident on the wall, $\bv\cdot\bn<0$, so that each ghost cell is known before the corresponding sweep begins.

The iteration is stopped when the $L^2$ norm of the residual of the steady equation \eqref{eq:steady_1d},
\begin{equation}
    \mathrm{Res}^{(\ell)}
    = \left(\int_0^1\!\int_{\bbR^3}\!\int_0^\infty
    \left| v_1\,\partial_x f^{(\ell)}
    - \frac{1}{\Kn}\,\mQ\big(f^{(\ell)},f^{(\ell)}\big)\right|^2
    dI\,d\bv\,dx\right)^{1/2},
    \label{eq:residual}
\end{equation}
falls below $10^{-5}$, where the spatial derivative is evaluated with the same upwind difference as in \eqref{eq:upwind}. In the two tests below we use $N_x=64$ cells, and the collision term is evaluated by the FLFSM with $N=72$ Fourier modes in each velocity direction and $K=10$ Laguerre modes, together with $S=16$, $L=(3+\sqrt2)S/4$, $\Lambda=S$, $N_\lambda=N$, and a Lebedev rule with $M=26$ angular nodes. Since the diffuse boundary conditions conserve mass but do not fix its value, we prescribe the total mass by $\int_0^1\rho\,dx=1$.

\subsubsection{Fourier flow}
\label{subsec:fourier_flow}

We first consider the classical Fourier flow, a standard benchmark for rarefied heat transfer in which the gas is confined between two stationary parallel plates held at different temperatures. The temperature difference drives a steady heat flux across the gap, and rarefaction produces temperature jumps at the walls. This problem has long been used to validate kinetic solvers \cite{yang1995rarefied,jaiswal2019discontinuous}.

In our test, both plates are stationary, and the walls at $x=0$ and $x=1$ are maintained at temperatures $\theta^W_L=0.8$ and $\theta^W_R=1.2$, respectively. We compute the steady solution for two Knudsen numbers, $\Kn=1$ and $\Kn=0.1$, and compare the density and temperature profiles against a DSMC reference computed with the same collision kernel.

\begin{figure}[!hptb]
\centering
\subfigure[$\rho$]{%
    \includegraphics[width=0.48\textwidth]{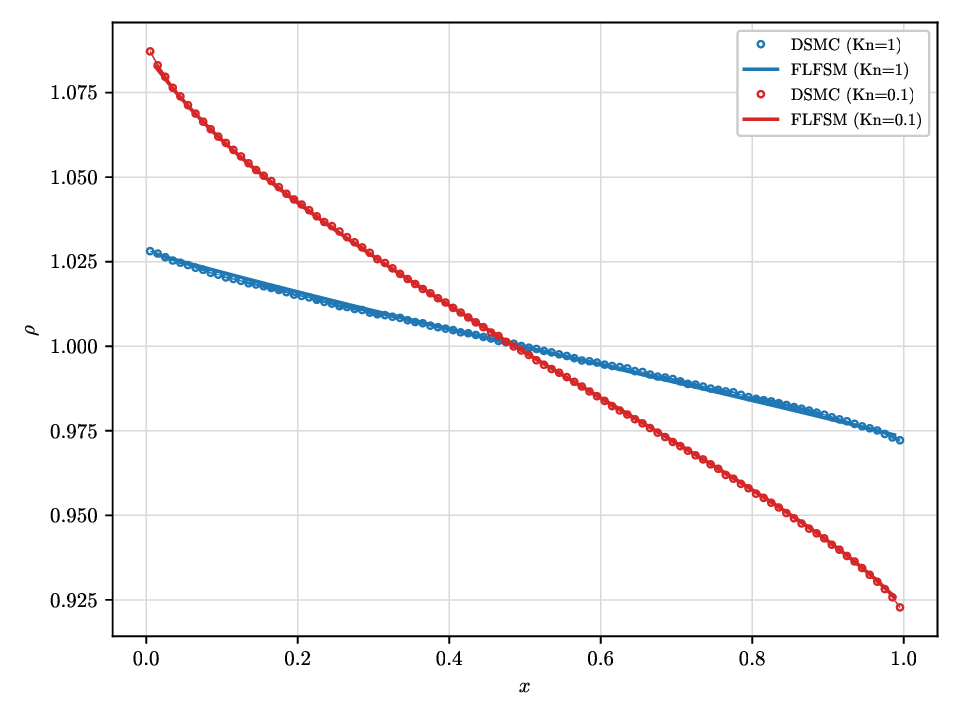}%
    \label{fig:fourier_rho}%
}
\subfigure[$\theta$]{%
    \includegraphics[width=0.48\textwidth]{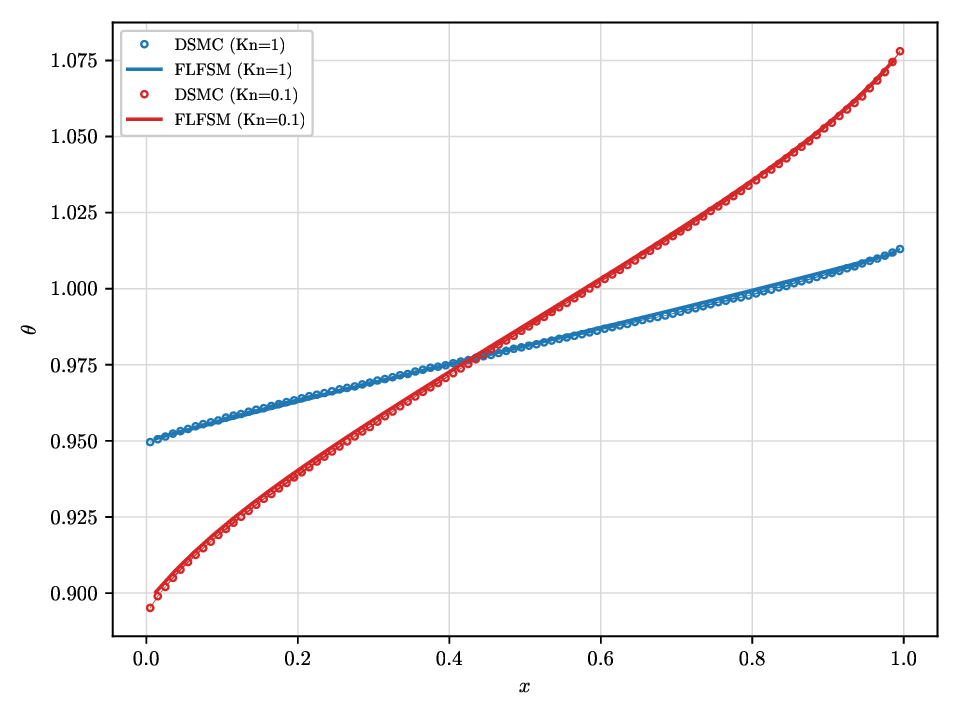}%
    \label{fig:fourier_T}%
}
\caption{(1D Fourier flow in Sec.~\ref{subsec:fourier_flow}) Steady-state numerical solutions of (a) the density $\rho$ and (b) the temperature $\theta$ for $\Kn=1$ and $0.1$. The solid lines are the FLFSM results and the circles are the DSMC solutions.}
\label{fig:fourier_flow}
\end{figure}

Fig.~\ref{fig:fourier_flow} shows the steady density $\rho$ and temperature $\theta$ obtained by the FLFSM, together with the DSMC results. For both $\Kn=1$ and $\Kn=0.1$, the FLFSM profiles agree well with DSMC, including the temperature jumps at the walls. The results confirm that the method captures the heat-transfer behavior of the polyatomic gas in different rarefaction regimes. The source iteration \eqref{eq:source_iteration} converges in $16$ iterations for $\Kn=1$ and in $46$ iterations for $\Kn=0.1$. As is typical of the source iteration, the number of iterations increases as $\Kn$ decreases.

\subsubsection{Couette flow}
\label{subsec:couette_flow}

We finally consider the planar Couette flow, which tests the FLFSM on a shear-driven flow. The two plates at $x=0$ and $x=1$ move tangentially in opposite directions with wall velocities $\bu^W=(0,-0.5,0)$ at $x=0$ and $\bu^W=(0,0.5,0)$ at $x=1$ and are both kept at the temperature $\theta^W=1$. Unlike the Fourier flow, the non-equilibrium here is driven by shear rather than heat conduction. A similar problem has been studied in \cite{shakhov1969couette,su2020can}. We again compute the steady solution for $\Kn=1$ and $\Kn=0.1$ and compare against a DSMC reference obtained with the same collision kernel.

\begin{figure}[!hptb]
\centering
\subfigure[$u_2$]{\includegraphics[width=0.45\textwidth]{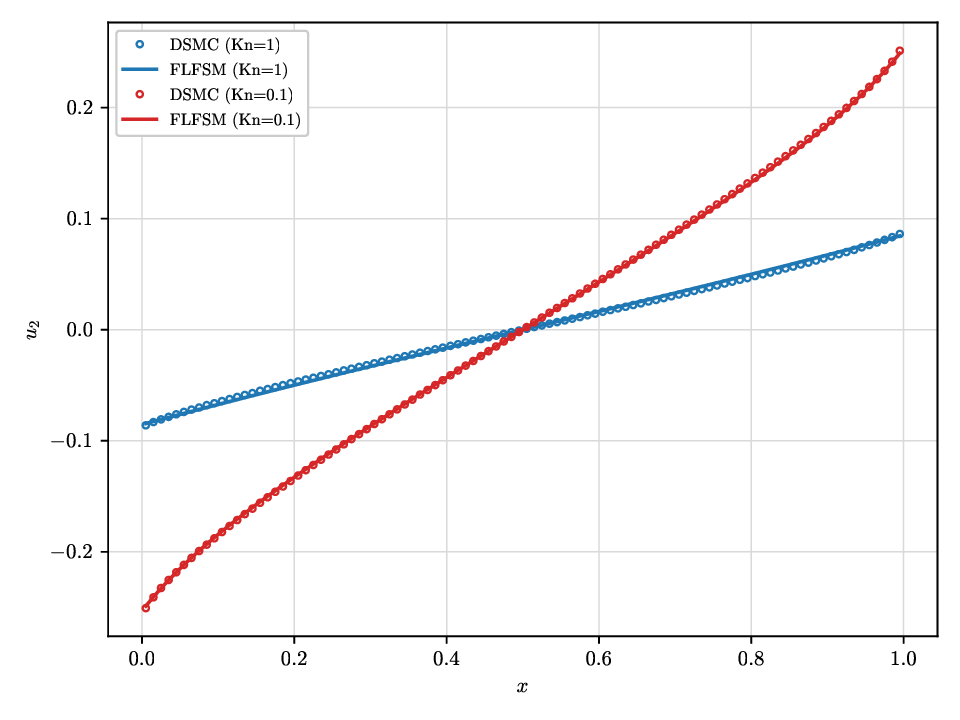}}
\subfigure[$q_1$]{\includegraphics[width=0.45\textwidth]{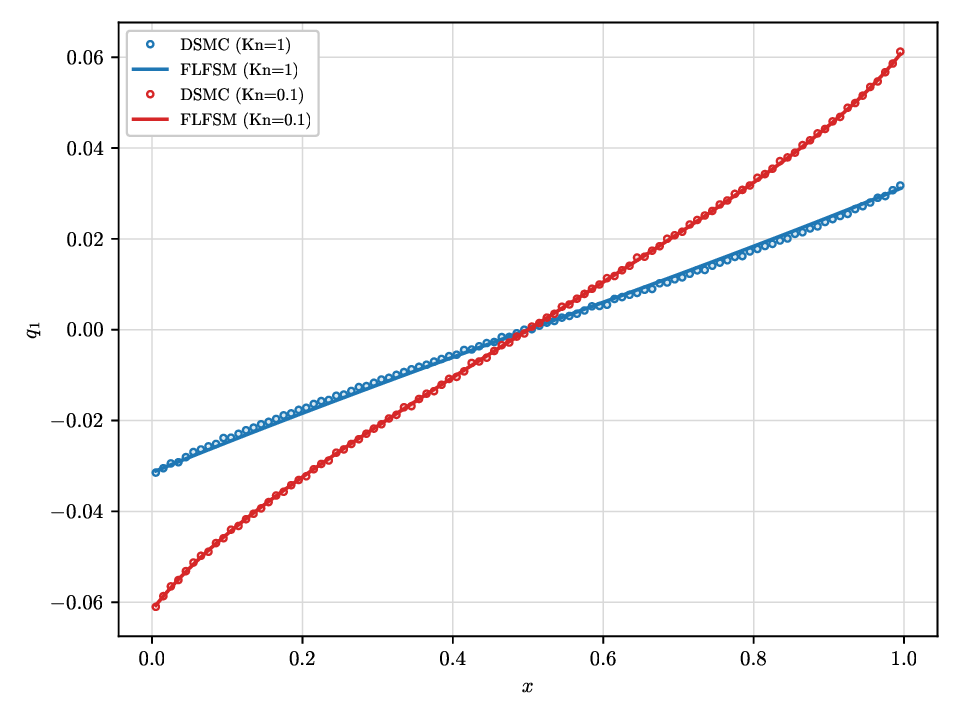}}
\caption{(1D Couette flow in Sec.~\ref{subsec:couette_flow}) Steady-state numerical solutions of the tangential velocity $u_2$ and the heat flux $q_1$ normal to the walls for $\Kn=1$ and $0.1$. The solid lines are the FLFSM results and the circles are the DSMC solutions.}
\label{fig:couette_flow}
\end{figure}

Fig.~\ref{fig:couette_flow} shows the steady tangential velocity $u_2$ and the heat flux $q_1$ obtained by the FLFSM, together with the DSMC results. For both $\Kn=1$ and $\Kn=0.1$, the FLFSM profiles agree well with DSMC. The source iteration converges in $8$ iterations for $\Kn=1$ and in $23$ iterations for $\Kn=0.1$, and, as in the Fourier flow, more iterations are needed for the smaller Knudsen number.

\section{Conclusion}
\label{sec:conclusion}

We have constructed a fast Laguerre--Fourier spectral method for the collision operator of the polyatomic Boltzmann equation with a continuous internal energy. The internal energy is discretized by generalized Laguerre functions whose weight matches that of the equilibrium, and the velocity by the Fourier spectral method. The essential step is a sequence of changes of variables that brings the gain term into a form in which the integrals over the internal energies and the energy-partition parameters are independent of the velocity convolutions. These integrals do not involve the distribution function and are therefore computed once and stored, which reduces the online cost of a single evaluation of the collision operator from $O(K^3N^6)$ for the direct spectral method to $O(K^2MN^4\log N+K^3MN^4)$, the same dependence on $N$ as in the monatomic fast spectral method.

The numerical experiments support this construction. The approximation of the collision operator converges spectrally in both the number of Fourier modes and the number of Laguerre modes. The measured online timings agree with \eqref{eq:cost_fast_vhs} on the finer velocity grids, while on the coarser ones they are dominated by the fixed overhead of each transform call. For the spatially homogeneous equation the FLFSM solutions agree with DSMC for an isotropic initial condition and two anisotropic ones, one of which is discontinuous. In addition, the two-temperature relaxation is validated for two collision kernels, against the exact macroscopic law for Maxwell molecules and against DSMC for hard spheres. Combined with a source iteration and a first-order upwind discretization in space, the method also agrees with DSMC on the one-dimensional Fourier and Couette flows at $\Kn=1$ and $\Kn=0.1$.

Finally, while this work focuses mainly on collision kernels depending only on the relative velocity, the energy-partition parameter $R$, and the collision angle, the extension of the method to the more general collision kernels in \cite{djordjic2023boltzmann}, which also depend on the energy-partition parameter $r$ and on the internal energies $I$ and $I_*$, is straightforward. Such an extension only changes the offline coefficients, while the online algorithm and its cost remain the same. Computations with more general collision kernels will produce results that better fit experiments. Numerical tests on these more realistic models are ongoing.

\appendix

\section{An exact two-temperature relaxation law}
\label{subsec:two_temperature_law}

We record here the moment calculation behind the relaxation reference used in
Sec.~\ref{subsec:relaxation_dsmc}. Consider the spatially homogeneous
equation
\begin{equation}
    \partial_t f = \mQ(f,f).
    \label{eq:homogeneous_poly}
\end{equation}
For a constant Maxwell-type kernel $\mathcal B=B_0$, normalize $\rho=1$ and
$\bu=\bz$. By the microreversibility of the polyatomic collision map
\cite{bourgat1994microreversible,djordjic2021consistent,gamba2023cauchy}, the
symmetric weak form is
\begin{equation}
\begin{aligned}
    \frac{d}{dt}\int_{\bbR^3}\int_0^\infty f\psi\,dI\,d\bv
    ={}& \frac12\int_{\bbR^3}\int_0^\infty\int_{\bbR^3}\int_0^\infty\int_{\bbS^2}\int_0^1\int_0^1 f f_*\,\mathcal B\,H_\delta \\
    &\times \big[\psi'+\psi_*'-\psi-\psi_*\big]\,dR\,dr\,d\boldsymbol{\sigma}\,dI_*\,d\bv_*\,dI\,d\bv ,
\end{aligned}
    \label{eq:weak_pair_form}
\end{equation}
where $\psi'=\psi(\bv',I')$ and $\psi_*'=\psi(\bv_*',I_*')$. Choosing the
collision invariants $\psi\in\{1,\bv,\tfrac12|\bv|^2+I\}$ gives
\begin{equation}
    \frac{d\rho}{dt}=0,\qquad \frac{d\bu}{dt}=\bz,\qquad \frac{d}{dt}\big(3\theta_{\mathrm{tr}}+\delta\theta_{\mathrm{int}}\big)=0.
    \label{eq:two_temperature_energy}
\end{equation}
The non-conserved part is obtained by taking $\psi=I$:
\begin{equation}
\begin{aligned}
    \frac{\delta}{2}\frac{d\theta_{\mathrm{int}}}{dt}
    ={}& \frac12\int_{\bbR^3}\int_0^\infty\int_{\bbR^3}\int_0^\infty\int_{\bbS^2}\int_0^1\int_0^1 f f_*\,\mathcal B\,H_\delta \\
    &\times \big[I'+I_*'-I-I_*\big]\,dR\,dr\,d\boldsymbol{\sigma}\,dI_*\,d\bv_*\,dI\,d\bv .
\end{aligned}
    \label{eq:internal_energy_weak}
\end{equation}

Using \eqref{eq:poly_post_internal_energy},
$I'+I_*'=(1-R)E$, this identity becomes
\begin{equation}
\begin{aligned}
    \frac{\delta}{2}\frac{d\theta_{\mathrm{int}}}{dt}
    ={}& \frac{B_0}{2}\int_{\bbR^3}\int_0^\infty\int_{\bbR^3}\int_0^\infty\int_{\bbS^2}\int_0^1\int_0^1 f f_*\,H_\delta \\
    &\times \big[(1-R)E-(I+I_*)\big]\,dR\,dr\,d\boldsymbol{\sigma}\,dI_*\,d\bv_*\,dI\,d\bv .
\end{aligned}
    \label{eq:internal_energy_substituted}
\end{equation}
Since $\mathcal B=B_0$ is constant, the integrals over $\boldsymbol{\sigma}$, $r$, and $R$ reduce to constants. Define
\begin{equation}
\begin{aligned}
    \nu_0 &= B_0\int_{\bbS^2}\!\int_0^1\!\int_0^1 H_\delta\,dR\,dr\,d\boldsymbol{\sigma} \\
    &= B_0\int_{\bbS^2}\!d\boldsymbol{\sigma}\int_0^1 [r(1-r)]^{\frac{\delta}{2}-1}\,dr\int_0^1 (1-R)^{\delta-1}R^{1/2}\,dR \\
    &= 4\pi B_0\,\mathrm{B}\!\big(\tfrac{\delta}{2},\tfrac{\delta}{2}\big)\,\mathrm{B}\!\big(\tfrac32,\delta\big),
\end{aligned}
    \label{eq:nu0}
\end{equation}
where $\mathrm{B}$ denotes the Beta function. The same calculation with the
extra factor $1-R$ gives
\begin{equation}
    \int_{\bbS^2}\!\int_0^1\!\int_0^1 (1-R)\,H_\delta\,dR\,dr\,d\boldsymbol{\sigma}
    =\frac{\nu_0}{B_0}\frac{\int_0^1 (1-R)^{\delta}R^{1/2}\,dR}{\int_0^1 (1-R)^{\delta-1}R^{1/2}\,dR}
    =\frac{\nu_0}{B_0}\frac{2\delta}{2\delta+3}.
    \label{eq:BLmoments}
\end{equation}
Substituting \eqref{eq:BLmoments} into \eqref{eq:internal_energy_substituted} yields
\begin{equation}
    \frac{\delta}{2}\frac{d\theta_{\mathrm{int}}}{dt} = \frac{\nu_0}{2}\int_{\bbR^3}\!\int_0^\infty\!\int_{\bbR^3}\!\int_0^\infty f f_*\left[\frac{2\delta}{2\delta+3}\,E-(I+I_*)\right] dI_*\,d\bv_*\,dI\,d\bv .
    \label{eq:internal_energy_reduced}
\end{equation}
Since $\rho=1$ and $\bu=\bz$, the definitions
\eqref{eq:poly_total_energy}--\eqref{eq:poly_temperatures} imply
\begin{equation}
    \frac{\delta}{2}\frac{d\theta_{\mathrm{int}}}{dt} = \frac{\nu_0}{2}\left[\frac{2\delta}{2\delta+3}\Big(\tfrac32\theta_{\mathrm{tr}}+\delta\theta_{\mathrm{int}}\Big)-\delta\theta_{\mathrm{int}}\right] = \frac{\nu_0}{2}\,\frac{3\delta}{2\delta+3}\,(\theta_{\mathrm{tr}}-\theta_{\mathrm{int}}),
    \label{eq:internal_pair_production}
\end{equation}
and therefore
\begin{equation}
    \frac{d\theta_{\mathrm{int}}}{dt} = \frac{3\nu_0}{2\delta+3}\,(\theta_{\mathrm{tr}}-\theta_{\mathrm{int}}).
    \label{eq:dtheta_int}
\end{equation}

Together with \eqref{eq:two_temperature_energy}, this gives the closed scalar
equation for the temperature gap,
\begin{equation}
    \frac{d}{dt}\big(\theta_{\mathrm{tr}}-\theta_{\mathrm{int}}\big)
    =-\left(1+\frac{\delta}{3}\right)\frac{d\theta_{\mathrm{int}}}{dt}
    =-\mu\,\big(\theta_{\mathrm{tr}}-\theta_{\mathrm{int}}\big),\qquad
    \mu = \frac{3+\delta}{2\delta+3}\,\nu_0 .
    \label{eq:mu}
\end{equation}
Hence
\begin{equation}
    \theta_{\mathrm{tr}}(t)-\theta_{\mathrm{int}}(t) = \big[\theta_{\mathrm{tr}}(0)-\theta_{\mathrm{int}}(0)\big]\,e^{-\mu t}, \qquad \theta_{\mathrm{tr/int}}(t) = \theta_\infty + \big[\theta_{\mathrm{tr/int}}(0)-\theta_\infty\big]\,e^{-\mu t},
    \label{eq:twoT_analytic}
\end{equation}
with
\begin{equation}
    (3+\delta)\,\theta_\infty = 3\theta_{\mathrm{tr}}(0)+\delta\theta_{\mathrm{int}}(0) .
    \label{eq:theta_infty}
\end{equation}

\section{Offline Evaluation of Internal Coefficients}
\label{app:offline_coefficients}

This appendix records the offline quadratures used to compute the coefficients in Sec.~\ref{subsec:fast_gain_loss}. The redistribution coefficient is evaluated for all needed values of $Y$ by
\begin{equation}
    G_{b,d}(Y)\approx \sum_{m=1}^{N_r} w_m^r\,
    \Phi_b(r_mY)\,\Phi_d((1-r_m)Y),
    \label{eq:Gbd_discrete}
\end{equation}
using a quadrature on $[0,1]$ adapted to the endpoint behavior of the trial functions. The coefficient $\mathcal P_a$ is evaluated by a Gauss--Jacobi rule for the weight $[\tau(1-\tau)]^{\delta/2-1}$,
\begin{equation}
    \mathcal P_a(U)\approx U^{\delta-1}\sum_{m=1}^{N_\tau} w_m^\tau\,
    \Psi_a(\tau_m U).
    \label{eq:Pa_discrete}
\end{equation}
For the kernel \eqref{eq:vhs_kernel}, the weight function \eqref{eq:F_weight} factorizes as
\begin{equation}
    F(g,R,\boldsymbol{\sigma},\boldsymbol{\xi}) = g^{2+\zeta}R^{\zeta/2}A(g,\boldsymbol{\sigma},\boldsymbol{\xi}), \qquad
    A(g,\boldsymbol{\sigma},\boldsymbol{\xi}) = \int_{\bbS^2}\mathfrak b(\boldsymbol{\omega}\cdot\boldsymbol{\sigma})\,e^{-\frac{i}{2}g\boldsymbol{\omega}\cdot\boldsymbol{\xi}}\,d\boldsymbol{\omega} ,
    \label{eq:F_factor}
\end{equation}
so that, with $R=\lambda^2/(\lambda^2+4Y)$, the coefficient \eqref{eq:K_abd} becomes
\begin{equation*}
    \mathcal{K}_{a,b,d}(g,\lambda,\boldsymbol{\sigma},\boldsymbol{\xi})=g^{2+\zeta}A(g,\boldsymbol{\sigma},\boldsymbol{\xi})\,\widetilde{\mathcal K}_{a,b,d}(g,\lambda),
\end{equation*}
where $\widetilde{\mathcal K}_{a,b,d}$ is independent of $\boldsymbol{\sigma}$ and $\boldsymbol{\xi}$. The semi-infinite $Y$-integral is evaluated without truncation by a Gauss--Laguerre rule after the shift $Y=Y_0(g_l,\lambda_q)+t$, $t\in[0,\infty)$,
\begin{equation}
\begin{aligned}
    \widetilde{\mathcal K}_{a,b,d}(g_l,\lambda_q) \approx{}&
    \sum_{n=1}^{N_Y} w_n^Y\,Y_n^{lq}
    \left(\frac{\lambda_q^2}{4}+Y_n^{lq}\right)^{-\delta-\frac12}
    \left(\frac{\lambda_q^2}{\lambda_q^2+4Y_n^{lq}}\right)^{\zeta/2} \\
    &\times \mathcal P_a\!\left(\frac{\lambda_q^2}{4}+Y_n^{lq}-\frac{g_l^2}{4}\right)
    G_{b,d}(Y_n^{lq}),
\end{aligned}
    \label{eq:Kabd_discrete}
\end{equation}
where $Y_n^{lq}=Y_0(g_l,\lambda_q)+t_n$ and $w_n^Y=\omega_n e^{t_n}$, with $\{t_n,\omega_n\}_{n=1}^{N_Y}$ the nodes and weights of the $N_Y$-point Gauss--Laguerre rule on $[0,\infty)$.

Once $\widetilde{\mathcal K}_{a,b,d}$ is available, the weight \eqref{eq:W_abd} is evaluated by the $g$-quadrature,
\begin{equation}
    W_{a,b,d}(\lambda_q,\boldsymbol{\sigma}_p,\boldsymbol{\xi}_{\boldsymbol{k}})
    \approx \frac{\lambda_q^2}{4}
    \sum_{l=1}^{N_g} w_l^g\,g_l^{2+\zeta}\,
    A(g_l,\boldsymbol{\sigma}_p,\boldsymbol{\xi}_{\boldsymbol{k}})\,
    \widetilde{\mathcal K}_{a,b,d}(g_l,\lambda_q).
    \label{eq:Wq}
\end{equation}
For the VHS-type kernel \eqref{eq:vhs_iso}, $A(g,\boldsymbol{\sigma},\boldsymbol{\xi})=4\pi B_0\operatorname{sinc}(g|\boldsymbol{\xi}|/2)$, so that $W_{a,b,d}$ is independent of $\boldsymbol{\sigma}_p$, depends on $\boldsymbol{\xi}_{\boldsymbol{k}}$ only through $|\boldsymbol{\xi}_{\boldsymbol{k}}|$, and is tabulated on the distinct radii of the Fourier grid.

The loss kernel \eqref{eq:ls_ker} is precomputed on the same grid. Since $\int_{\bbS^2}\mathfrak b(\hat{\bg}\cdot\boldsymbol{\sigma})\,d\boldsymbol{\sigma}=2\pi\int_{-1}^1\mathfrak b(s)\,ds$ is independent of $\hat{\bg}$, integrating over the direction of $\bg$ yields a sinc factor, and \eqref{eq:ls_ker} reduces to the radial quadrature
\begin{equation}
    \widehat{\mathcal B}_{\boldsymbol{k}}
    \approx 8\pi^2\,c_\delta\int_{-1}^1\mathfrak b(s)\,ds \sum_{l=1}^{N_g} w_l^g\,g_l^{2+\zeta}\,
    \operatorname{sinc}\!\big(g_l|\boldsymbol{\xi}_{\boldsymbol{k}}|\big),
    \label{eq:Bhat_discrete}
\end{equation}
where $c_\delta=\int_0^1\!\int_0^1 R^{\zeta/2}H_\delta(r,R)\,dr\,dR=\mathrm{B}\big(\tfrac{\delta}{2},\tfrac{\delta}{2}\big)\,\mathrm{B}\big(\tfrac{\zeta+3}{2},\delta\big)$, and $\mathrm B$ denotes the Beta function.

\bibliographystyle{siamplain}
\bibliography{article}

@book{bird1994molecular,
	series = {Oxford {Engineering} {Science} {Series}},
	title = {Molecular Gas Dynamics and the Direct Simulation of Gas Flows},
	publisher = {Oxford University Press},
	author = {Bird, G.},
	year = {1994}
}

@article{gamba2017fast,
	title = {A fast spectral method for the {Boltzmann} collision operator with general collision kernels},
	volume = {39},
	number = {4},
	journal = {SIAM J. Sci. Comput.},
	author = {Gamba, I. and Haack, J. and Hauck, C. and Hu, J.},
	year = {2017},
	pages = {B658--B674}
}

@article{mouhot2006fast,
	title = {Fast algorithms for computing the {Boltzmann} collision operator},
	volume = {75},
	number = {256},
	journal = {Math. Comput.},
	author = {Mouhot, C. and Pareschi, L.},
	year = {2006},
	pages = {1833-1852}
}

@article{wu2013deterministic,
	title = {Deterministic numerical solutions of the {Boltzmann} equation using the fast spectral method},
	volume = {250},
	journal = {J. Comput. Phys.},
	author = {Wu, L. and White, C. and Scanlon, T. and Reese, J. and Zhang, Y.},
	year = {2013},
	pages = {27--52}
}

@article{wang2019approximation,
	title = {Approximation of the {Boltzmann} collision operator based on {Hermite} spectral method},
	journal = {J. Comput. Phys.},
	author = {Wang, Y. and Cai, Z.},
	year = {2019},
        volume = {397},
	pages = {108815}
}

@book{shen2006rarefied,
  title={Rarefied Gas Dynamics: Fundamentals, Simulations and Micro Flows},
  author={Shen, C.},
  year={2006},
  publisher={Springer Science \& Business Media}
}

@article{broadwell1964study,
  title={Study of rarefied shear flow by the discrete velocity method},
  author={Broadwell, J. E.},
  journal={J. Fluid Mech.},
  volume={19},
  number={3},
  pages={401--414},
  year={1964},
  publisher={Cambridge University Press}
}

@article{buet1996discrete,
  title={{A discrete-velocity scheme for the Boltzmann operator of rarefied gas dynamics}},
  author={Buet, C.},
  journal={Transp. Theory Stat. Phys.},
  volume={25},
  number={1},
  pages={33--60},
  year={1996},
  publisher={Taylor \& Francis}
}

@article{li2023hermite,
  title={Hermite spectral method for the inelastic {Boltzmann} equation},
  author={Li, R. and Lu, Y. and Wang, Y.},
  journal={Phys. Fluids},
  volume={35},
  number={10},
  year={2023},
  publisher={AIP Publishing}
}

@article{yang1995rarefied,
  title={Rarefied flow computations using nonlinear model {Boltzmann} equations},
  author={Yang, J. and Huang, J.},
  journal={J. Comput. Phys.},
  volume={120},
  number={2},
  pages={323--339},
  year={1995},
  publisher={Elsevier}
}

@article{maxwell1878iii,
  title={On stresses in rarefied gases arising from inequalities of temperature},
  author={Maxwell, J.C.},
  journal={Proc. R. Soc. Lond.},
  volume={27},
  number={185-189},
  pages={304--308},
  year={1878},
  publisher={The Royal Society London}
}

@article{shakhov1969couette,
  title={Couette problem for the generalized {Krook} equation stress-peak effect},
  author={Shakhov, E.M.},
  journal={Fluid Dyn.},
  volume={4},
  number={5},
  pages={9--13},
  year={1969},
  publisher={Springer}
}

@article{hu2020burnett,
  title={Burnett spectral method for high-speed rarefied gas flows},
  author={Hu, Z. and Cai, Z.},
  journal={SIAM J. Sci. Comput.},
  volume={42},
  number={5},
  pages={B1193--B1226},
  year={2020},
  publisher={SIAM}
}

@article{goldstein1989investigations,
  title={Investigations of the motion of discrete-velocity gases},
  author={Goldstein, D. and Sturtevant, B. and Broadwell, J. E.},
  journal={Prog. Astronaut. Aeronaut.},
  volume={117},
  pages={100--117},
  year={1989}
}

@article{cai2020burnett,
  title={Burnett spectral method for the spatially homogeneous {Boltzmann} equation},
  author={Cai, Z. and Fan, Y. and Wang, Y.},
  journal={Comput. Fluids},
  volume={200},
  pages={104456},
  year={2020},
  publisher={Elsevier}
}

@article{hu2022fast,
  title={A Fast {Petrov--Galerkin} Spectral Method for the Multidimensional {Boltzmann} Equation Using Mapped {C}hebyshev Functions},
  author={Hu, J. and Huang, X. and Shen, J. and Yang, H.},
  journal={SIAM J. Sci. Comput.},
  volume={44},
  number={3},
  pages={A1497--A1524},
  year={2022},
  publisher={SIAM}
}

@article{hu2020petrov,
  title={A {Petrov-Galerkin} spectral method for the inelastic {Boltzmann} equation using mapped {Chebyshev} functions},
  author={Hu, J. and Shen, J. and Wang, Y.},
  journal={Kinet. Relat. Models},
  volume={13},
  number={4},
  pages = {677--702},
  year={2020}
}

@article{pareschi2000numerical,
  title={Numerical solution of the {Boltzmann} equation {I}: Spectrally accurate approximation of the collision operator},
  author={Pareschi, L. and Russo, G.},
  journal={SIAM J. Numer. Anal.},
  volume={37},
  number={4},
  pages={1217--1245},
  year={2000},
  publisher={SIAM}
}

@article{pareschi1996fourier,
  title={A {Fourier} spectral method for homogeneous {Boltzmann} equations},
  author={Pareschi, L. and Perthame, B.},
  journal={Transp. Theory Stat. Phys.},
  volume={25},
  number={3-5},
  pages={369--382},
  year={1996},
  publisher={Taylor \& Francis}
}

@article{yin2025fast,
author = {Yin, T. and Cai, Z. and Wang, Y.},
title = {A Fast {Fourier} Spectral Method for the Linearized {Boltzmann} Collision Operator},
 journal={SIAM J. Sci. Comput.},
volume = {48},
number = {1},
pages = {B83-B111},
year = {2026},
doi = {10.1137/25M1737791}
}

@article{su2020can,
  title={Can we find steady-state solutions to multiscale rarefied gas flows within dozens of iterations?},
  author={Su, W. and Zhu, L. and Wang, P. and Zhang, Y. and Wu, L.},
  journal={J. Comput. Phys.},
  volume={407},
  pages={109245},
  year={2020},
  publisher={Elsevier}
}

@article{jaiswal2019multispecies,
  title={A discontinuous {Galerkin} fast spectral method for the multi-species {Boltzmann} equation},
  author={Jaiswal, S. and Alexeenko, A. A. and Hu, J.},
  journal={Comput. Methods Appl. Mech. Engrg.},
  volume={352},
  pages={56--84},
  year={2019},
  doi={10.1016/j.cma.2019.04.015}
}

@article{jaiswal2019discontinuous,
  title={A discontinuous {Galerkin} fast spectral method for the full {Boltzmann} equation with general collision kernels},
  author={Jaiswal, S. and Alexeenko, A. and Hu, J.},
  journal={J. Comput. Phys.},
  volume={378},
  pages={178--208},
  year={2019},
  publisher={Elsevier}
}

@article{filbet2006solving,
  title={Solving the {Boltzmann} equation in $N\log_2 N$},
  author={Filbet, F. and Mouhot, C. and Pareschi, L.},
  journal={SIAM J. Sci. Comput.},
  volume={28},
  number={3},
  pages={1029--1053},
  year={2006},
  publisher={SIAM}
}

@article{pareschi2000fokker,
  title={Fast spectral methods for the {Fokker}--{Planck}--{Landau} collision operator},
  author={Pareschi, L. and Russo, G. and Toscani, G.},
  journal={J. Comput. Phys.},
  volume={165},
  number={1},
  pages={216--236},
  year={2000},
  publisher={Elsevier}
}

@article{hu2012fast,
  title={A fast spectral algorithm for the quantum {Boltzmann} collision operator},
  author={Hu, J. and Ying, L.},
  journal={Commun. Math. Sci.},
  volume={10},
  number={3},
  pages={989--999},
  year={2012},
  publisher={International Press}
}

@article{hu2019fast,
  title={A fast spectral method for the inelastic {Boltzmann} collision operator and application to heated granular gases},
  author={Hu, J. and Ma, Z.},
  journal={J. Comput. Phys.},
  volume={385},
  pages={119--134},
  year={2019},
  publisher={Elsevier}
}

@article{borgnakke1975statistical,
  title   = {Statistical collision model for {Monte} {Carlo} simulation of polyatomic gas mixture},
  author  = {Borgnakke, C. and Larsen, P. S.},
  journal = {J. Comput. Phys.},
  volume  = {18},
  number  = {4},
  pages   = {405--420},
  year    = {1975},
  doi     = {10.1016/0021-9991(75)90094-7}
}

@article{bourgat1994microreversible,
  title   = {Microreversible collisions for polyatomic gases and {Boltzmann's} theorem},
  author  = {Bourgat, J.-F. and Desvillettes, L. and Le Tallec, P. and Perthame, B.},
  journal = {Eur. J. Mech. B Fluids},
  volume  = {13},
  number  = {2},
  pages   = {237--254},
  year    = {1994}
}

@article{djordjic2021consistent,
  title   = {Consistent, explicit, and accessible {Boltzmann} collision operator for polyatomic gases},
  author  = {Djordji\'c, V. and Pavi\'c-\v{C}oli\'c, M. and Torrilhon, M.},
  journal = {Phys. Rev. E},
  volume  = {104},
  number  = {2},
  pages   = {025309},
  year    = {2021},
  doi     = {10.1103/PhysRevE.104.025309}
}

@article{djordjic2023boltzmann,
  title   = {{Boltzmann} collision operator for polyatomic gases in agreement with experimental data and {DSMC} method},
  author  = {Djordji\'c, V. and Oblapenko, G. and Pavi\'c-\v{C}oli\'c, M. and Torrilhon, M.},
  journal = {Contin. Mech. Thermodyn.},
  volume  = {35},
  number  = {1},
  pages   = {103--119},
  year    = {2023},
  doi     = {10.1007/s00161-022-01167-8}
}

@article{gamba2023cauchy,
  title   = {On the {Cauchy} problem for {Boltzmann} equation modeling a polyatomic gas},
  author  = {Gamba, I. M. and Pavi\'c-\v{C}oli\'c, M.},
  journal = {J. Math. Phys.},
  volume  = {64},
  number  = {1},
  pages   = {013303},
  year    = {2023},
  doi     = {10.1063/5.0103621}
}

@article{munafo2014spectral,
  title   = {A spectral-{Lagrangian} {Boltzmann} solver for a multi-energy level gas},
  author  = {Munaf\`o, A. and Haack, J. R. and Gamba, I. M. and Magin, T. E.},
  journal = {J. Comput. Phys.},
  volume  = {264},
  pages   = {152--176},
  year    = {2014},
  doi     = {10.1016/j.jcp.2014.01.009}
}

@book{shen2011spectral,
  title={Spectral Methods: Algorithms, Analysis and Applications},
  author={Shen, J. and Tang, T. and Wang, L.-L.},
  series={Springer Series in Computational Mathematics},
  volume={41},
  year={2011},
  publisher={Springer, Berlin}
}

@article{rykov1975model,
  title   = {A model kinetic equation for a gas with rotational degrees of freedom},
  author  = {Rykov, V. A.},
  journal = {Fluid Dyn.},
  volume  = {10},
  number  = {6},
  pages   = {959--966},
  year    = {1975},
  doi     = {10.1007/BF01023275}
}

@article{mathiaud2022esbgk,
  title   = {An {ES-BGK} model for diatomic gases with correct
             relaxation rates for internal energies},
  author  = {Mathiaud, J. and Mieussens, L. and Pfeiffer, M.},
  journal = {Eur. J. Mech. B Fluids},
  volume  = {96},
  pages   = {65--77},
  year    = {2022},
  doi     = {10.1016/j.euromechflu.2022.07.003}
}

@article{arima2024novel,
  title   = {A novel {ES-BGK} model for non-polytropic gases
             with internal state density independent of the temperature},
  author  = {Arima, T. and Mentrelli, A. and Ruggeri, T.},
  journal = {J. Stat. Phys.},
  volume  = {191},
  pages   = {95},
  year    = {2024},
  doi     = {10.1007/s10955-024-03286-9}
}

@article{zeng2023general,
  title   = {General synthetic iterative scheme for polyatomic
             rarefied gas flows},
  author  = {Zeng, J. and Yuan, R. and Zhang, Y. and Li, Q. and Wu, L.},
  journal = {Comput. Fluids},
  volume  = {265},
  pages   = {105998},
  year    = {2023},
  doi     = {10.1016/j.compfluid.2023.105998}
}

\end{document}